\documentclass[10pt,journal,twoside]{IEEEtran}
\usepackage{booktabs} 
\usepackage{algorithm2e}
\usepackage{amsfonts,amssymb}
\usepackage{hyperref}

\usepackage[caption=false]{subfig}
\usepackage{multirow, makecell}
\usepackage[dvipsnames]{xcolor}
\usepackage[prependcaption,textsize=small]{todonotes}
\usepackage{diagbox}

\DeclareMathAlphabet{\mathpzc}{OT1}{pzc}{m}{it}

\newcommand{\z}{\boldsymbol{z}}
\newcommand{\x}{\boldsymbol{x}}

\newcommand{\vbold}{\boldsymbol{v}}
\newcommand{\ubold}{\boldsymbol{u}}

\newcommand{\psibold}{\boldsymbol{\psi}}
\newcommand{\y}{\boldsymbol{y}}

\newcommand{\Q}{\mathcal Q}

\newcommand{\G}{\mathcal G}
\newcommand{\N}{\mathcal N}
\newcommand{\E}{\mathcal E}
\newcommand{\B}{\mathcal B}
\newcommand{\K}{\mathcal K}
\newcommand{\U}{\mathcal U}
\newcommand{\T}{\mathcal T}

\newcommand{\A}{\mathcal{A}}
\newcommand{\V}{\mathcal{V}}
\newcommand{\underl}{\underline}
\newcommand{\overl}{\overline}

\usepackage{cite}

\ifCLASSINFOpdf
\else
\fi
\usepackage{amsmath}
\newtheorem{theorem}{Theorem}
\newtheorem{remark}{Remark}

\begin{document}
%
\title{Unit Commitment with Gas Network Awareness}
%
%
%

\author{Geunyeong Byeon~
        and Pascal Van Hentenryck,~\IEEEmembership{Member,~IEEE,}
\thanks{G. Byeon is with the Department
of Industrial and Operations Engineering, University of Michigan, Ann Arbor,
MI, 48109 (e-mail: gbyeon@umich.edu).}
\thanks{P. Van Hentenryck is with the H. Milton Stewart School of Industrial and Systems Engineering, Georgia Institute of Technology, Atlanta, GA, 30332 (e-mail: pvh@isye.gatech.edu).}
}

%
%

\markboth{Byeon and Van Hentenryck: Unit Commitment with Gas Network Awareness}%
{Byeon and Van Hentenryck: Unit Commitment with Gas Network Awareness}
%



\maketitle

\begin{abstract}
  Recent changes in the fuel mix for electricity generation and, in
  particular, the increase in Gas-Fueled Power Plants (GFPP), have
  created significant interdependencies between the electrical power
  and natural gas transmission systems. However, despite their
  physical and economic couplings, these networks are still operated
  independently, with asynchronous market mechanisms. This mode of
  operation may lead to significant economic and reliability risks in
  congested environments as revealed by the 2014 polar vortex event
  experienced by the northeastern United States. To mitigate these
  risks, while preserving the current structure of the markets, this
  paper explores the idea of introducing gas network awareness into
  the standard unit commitment model. Under the assumption that the
  power system operator has some (or full) knowledge of gas demand
  forecast and the gas network, the paper proposes a tri-level
  mathematical program where natural gas zonal prices are given by the
  dual solutions of natural-gas flux conservation constraints and
  commitment decisions are subject to bid-validity constraints that
  ensure the economic viability of the committed GFPPs. This tri-level
  program can be reformulated as a single-level Mixed-Integer
  Second-Order Cone program which can then be solved using a dedicated
  Benders decomposition. The approach is validated on a case study for
  the Northeastern United States \cite{bent2018joint} that can
  reproduce the gas and electricity price spikes experienced during
  the early winter of 2014. The results on the case study demonstrate
  that gas awareness in unit commitment is instrumental in avoiding
  the peaks in electricity prices while keeping the gas prices to
  reasonable levels.
  \end{abstract}


%
\IEEEpeerreviewmaketitle

\section{Introduction}\label{sec:Intro}
%
%
%
%

\IEEEPARstart{G}{as}-Fueled Power Plants (GFPPs) have become a
significant part of the energy mix in the last decades, primarily
because of their operational flexibility and lower environmental
impacts. Although GFPPs have introduced interdependencies between the
natural gas and electrical power systems, these networks are still
operated independently, with asynchronous market mechanisms. In
particular, the unit commitment decisions in the electrical power
system take place before the realization of natural gas spot prices,
introducing reliability risks and economic inefficiencies in congested
environments. Indeed, the GFPPs may not be able to secure gas at
reasonable prices, introducing either reliability issues or
electricity gas spikes.

This undesirable outcome occurred in the Northeastern United States
during the early winter of 2014. Extremely low temperatures induced an
unusual coincident peak in electricity and natural gas demand. On the
one hand, it produced record-high natural gas spot prices due to
congestion. On the other hand, high electricity loads led the
electrical power system operator to call for some emergency actions,
which resulted in higher electricity prices
\cite{pgm2014analysis}. Moreover, the power system operator, valuing
reliability the most, encouraged committed GFPPs to buy natural gas at
all costs without assurance of cost recovery, further aggravating the
economic cost \cite{FERC2015}. It is important to mention that the
critical issue in this case was not the gas supply, but rather congestion
in the gas transmission network. Moreover, a recent study
\cite{bent2018joint} has shown that the cost of expanding the gas and
network infrastructures to avoid such events would be prohibitive.



To address these interdependencies, a number of researchers have
studied how to incorporate the natural gas transmission capabilities
into the operational decisions of electrical power systems. See, for
instance,
\cite{liu2009security,liu2011coordinated,martinez2012unified,correa2014security,correa2015integrated,
  alabdulwahab2015coordination, biskas2016coupled, li2017security,
  zhao2017unit, zlotnik2017coordinated}. Other researchers have also
studied how to incorporate the economic coupling between these two
infrastructures using new market mechanisms. A new market framework
with a joint ISO, using price- or volume-based approaches, was
investigated in
\cite{chen2018clearing,ordoudis2017exploiting}. Instead of introducing
one joint ISO, other researchers have proposed a new market framework
that assumes centralized independent gas markets, synchronizes the
electricity and gas market days, and allows some information exchange
between some parties in the electricity and gas markets (e.g., market
operators or GFPPs)
\cite{gil2016electricity,chen2017multi,wang2018equilibrium,ji2017coordinated,ji2018day,zhao2018shadow}.

This paper takes a different approach that stays within the current
operating practices and does not introduce a new market
mechanism. Instead, the approach generalizes the unit commitment model
to capture the physical and economic couplings and strive to ensure
both physical feasibility and economic viability. More precisely, the
paper introduces the Unit Commitment problem with Gas Network
Awareness (UCGNA) to schedule a set of generating units for the next
day while taking account the fuel delivery and the natural gas
prices that are propagated back by the natural gas system. The UCGNA
imposes bid-validity constraints on the GFPPs to ensure their
profitability and estimates the natural gas prices for these
constraints with the dual solutions associated with the flux
conservation constraints of the gas market.

The UCGNA is formulated as a tri-level mathematical program and
assumes that the power system operator has partial (or full) knowledge
on gas demand forecast and gas network. When the power system is
modeled with its DC approximation and the gas network with the
second-order cone program from \cite{sanchez2016convex} to model its
steady-state physics, the tri-level mathematical program can be
reformulated as a single-level Mixed-Integer Second-Order Cone Program
(MISOCP) through strong duality of the innermost problem. The
resulting MISOCP can then be solved using a dedicated Benders
decomposition recently proposed in \cite{ByeonTheory2018}.

The key contributions of this paper are threefold. First, it proposes
the first unit commitment model (UCGNA) that incorporates both the
physical and economic couplings of electrical power and natural gas
transmission systems and can be used within current operating
practices. Second, it proposes a MISOCP that captures the UCGNA and
can be solved through Benders decomposition. Finally, it demonstrates
the potential of the approach on a detailed case study that replicates
the behavior of the 2014 polar vortex event on the Northeastern United
States. In particular, the paper shows that, on the case study, the
UCGNA avoids the electricity price peaks and keeping the total gas
costs reasonable, contrary to current practice, even for highly congested
electrical and gas networks. 

The rest of this paper is organized as follows. Section
\ref{sec:UCGNA} formalizes the UCGNA and Section \ref{sec:approxi}
presents the MISOCP. Section \ref{sec:sol_approach} briefly reviews
the solution methods for the MISOCP. Section \ref{sec:data} describes
the test cases. Lastly, Section \ref{sec:case_study} analyzes the
behavior of the model on the case study and Section
\ref{sec:conclusion} concludes the paper.

\section{Unit Commitment With Gas Awareness}
\label{sec:UCGNA}


This section specifies the UCGNA, including its electricity system,
its natural gas network, and their physical and economic
couplings. The electricity transmission grid is represented by an
undirected graph $ \G^e = (\N, \E)$ and the natural gas transmission
system by a directed graph $\G^g = (\V, \A)$. Boldface letters
represent vectors of variables, $[a,b]_{\mathbb Z}$ denotes the set of
integers in interval $[a,b]$, and $[n]$ denotes the set
$\{1,\cdots,n\}$ for some integer $n \ge 1$. The letter $\T$ denotes
the set of time periods $\{0,1, \cdots, T\}$.



\subsection{The Electricity Transmission System}
\label{sec:e}

In the United States, Unit Commitment (UC) and Economic Dispatch (ED)
problems are solved daily to determine the hourly operating schedule
of generating units for the next day from bids submitted by market
participants. Tables \ref{table:param:e} and \ref{table:var:e}
summarize the parameters and variables of the UC/ED problems. With
these notations, the UC model is specified in Figure \ref{fig:UCED}:
It is standard but is presented as a bi-level program to make the
UCGNA formulations more intuitive subsequently.

\begin{table}[!t]\fontsize{9}{9}\selectfont
	\centering
	\caption{Parameters of the Electricity System.} \label{table:param:e}
	{\begin{tabular}{p{0.1\textwidth} p{0.35\textwidth}}
			\toprule
			$\G^e = (\N, \E)$  & Undirected graph where $\N$ is a set of buses indexed by $i = 1, \cdots, N$ and $\E$ is a set of lines indexed with $e = 1, \cdots, E$ \\%
			$\U$ & Set of generators, indexed by $u = 1, \cdots, U$\\%
			$\quad\U^g \subseteq \U$ & Set of GFPPs \\
			$\quad\U (i) \subseteq \U$ & Set of generators located at $i \in \N$ \\%
			$\B_u$ & Set of supply bids submitted by $u \in \U$, indexed by $b = 1, \cdots, B_u$ \\
			$\quad \beta_b$ & Bid price of $b \in \B_u$\\ 
			$\quad\overline {s}_{b}$& Amount of 
			real power generation of $b \in \B_u$ \\%
			$\underline p_u, \overline p_u$ &  Minimum/maximum real power generation of $u \in \U$ \\%
			$\underline R_u, \overline R_u$ & Ramp-down/-up rate of $u \in \U$ \\
			$c_u$ & No-load cost of $u \in \U$\\%
			$\Psi_{u}$ & Set of counts of time periods with distinct start-up costs of $u$ indexed by $h$ \\
			$\quad {C}_{u,h}$ & Start-up cost of $u\in \U$ when $u$ is turned on after it has been offline for some time $\in [\Psi_{u,h}, \Psi_{u,h+1}]$\\%
			$\overline u_{u,0}, \overline p_{u,0}$ & Initial on-off status/real power generation of $u \in \U$ \\%
			$ \underline\tau_{u}, \overline \tau_{u}$ &  Minimum-down/-up time of  $u \in \U$ \\  
			$\underline \tau_{u,0}, \overline \tau_{u,0}$ &  The time that generator $u \in \U$ has to be inactive/active from $t = 0$ \\%
			$b_{e}$ & Line susceptance of $e \in \E$\\%
			$\overline f_{e}$ & Real power limit of $e \in \E$\\%
			$(d^e_{i,t})_{i \in \N}$ & Electricity load profile during $t \in \T$ \\%
			$\Delta_e$ & Maximum voltage angle difference between two end-points of $e \in \E$\\
			$\underline \theta_i, \overline \theta_i$ & Minimum/maximum voltage angle at $i \in \N$\\
			\bottomrule
	\end{tabular}}{}
\end{table}

\begin{table}[!t]\fontsize{9}{9}\selectfont
	\centering
	\caption{Variables of the Electricity System.}
        \label{table:var:e}
	{\begin{tabular}{p{0.03\textwidth} p{0.42\textwidth}}
			\toprule
			\multicolumn{2}{l}{  \textbf{Binary variables}} \\
			$u_{u,t}$ &  1 if $u \in \U$ is on during $t \in \T$, 0 otherwise \\
			$v^+_{u,t}$ &  1 if $u \in \U$ becomes online during $t \in \T$, 0 otherwise \\
			$v^-_{u,t}$ &  1 if $u \in \U$ becomes offline during $t \in \T$, 0 otherwise \\
			%
			\multicolumn{2}{l}{ \textbf{Continuous variables}} \\
			$s^e_{b,t}$ & Real power generation from $b \in \B_u$ of $u \in \U$ during $t \in \T$\\
			$p_{u,t} $ & Real power generation of $u \in \U$ during $t \in \T$ \\
			$f_{e,t}$ & Real power flow on $e \in \E$ during $t \in \T$ \\
			$r_{u,t}$ & Start-up cost of $u \in \U$ during $t \in \T$ \\
			$\theta_{i,t}$ & Voltage angle on $i \in \N$ during $t \in \T$\\
			\bottomrule
	\end{tabular}}{}\end{table}

\begin{figure}[!t]
\begin{subequations}\fontsize{9}{9}\selectfont
		\begin{align}
		 \min \ &  \sum_{t\in [T]} \sum_{u \in \U}  (
		c_u u_{u,t} + r_{u,t} +\sum_{b \in \B_u} \beta_b {s}^e_{b,t} ) \label{e:obj}\\
		\mbox{s.t.} \  
		&r_{u,t} \ge C_{u,h} ( u_{u,t} - \sum_{n \in [h]} u_{u,t-n} ), \nonumber\\
		& \hspace{0.2\textwidth}  \forall h \in \Psi_s, u \in \U, t \in [T], \label{e:UC:su}\\
		&r_{u,t} \ge 0, \ \forall u \in \U, t \in [T], \label{e:UC:su:nonnegative}\\
		&u_{u,t} = \overline u_{u,0}, \ \forall u \in \U, \  t \in [0,\overline\tau_{u,0} + \underline \tau_{u,0}]_{\mathbb Z}, \label{e:UC:initial_gen_status}\\
		& \sum_{t' \in [t - \overline\tau_{u} +1, t]_{\mathbb Z}} v^+_{u,t'} \le u_{u,t},  \nonumber\\
		& \hspace{0.1\textwidth}  \forall  u \in \U, \ t \in [\max\{\overline\tau_{u}, \overline\tau_{u,0} + 1\}, T]_\mathbb{Z}, \label{e:UC:min_up}\\
		& \sum_{t' \in [t - \underline \tau_{u} +1,t]_{\mathbb Z} } v^+_{u,t'} \le 1- u_{u,t - \underline \tau_{u}},  \nonumber\\
		& \hspace{0.1\textwidth}  \forall u \in \U, \ t \in [\max\{\underline \tau_{u}, \underline \tau_{u,0} + 1\}, T]_\mathbb{Z}, \ \label{e:UC:min_down}\\
		& v^+_{u,t} - v^-_{u,t} = u_{u,t} - u_{u,t-1},\forall u \in \U, \  t \in [T], \label{e:UC:logic:on_off}\\
		& v^+_{u,t}, v^-_{u,t}, u_{u,t} \in \{0,1\},  \forall u \in \U, \ t \in [T], \label{e:UC:logic:nonnega}\\
		& \boldsymbol s^e = \mbox{argmin} \ \Q (\ubold, \vbold^+, \vbold^-), \label{eq:ED}
		\end{align}

where $\Q (\ubold, \vbold^+, \vbold^-)$ denotes the ED problem specified as follows:

		\begin{align}
		\min &  \sum_{t \in [T]} \sum_{u \in \U} \left(\sum_{b \in \B_u} \beta_b {s}^e_{b,t} \right) \label{e:ED:obj}\\
		  \mbox{s.t.} & \sum_{u \in \U(i)} {p}_{u,t} - d^e_{i,t}= \sum_{e \in \E: e_t = i} f_{e,t} - \sum_{e \in \E: e_h = i} f_{e,t}, \nonumber \\
		  &\hspace{0.25\textwidth} \ \forall i \in \N,\ t \in [T], \label{e:ED:bal}\\
		  & {p}_{u,t} = \sum_{b \in \B_u} {s}^e_{b,t}  \ \forall u \in \U, \ t \in [T], \label{e:ED:bid}\\
		  &  0 \le s^e_{b,t} \le \overline {s}_{b},  \ \forall b \in \B_u, \ u \in \U, \ t \in [T], \label{e:ED:bound:bids} \\
		  &  \underline{p}_{u}u_{u,t} \le p_{u,t} \le \overline{p}_{u} u_{u,t},  \ \forall u \in \U, \ t \in [T], \label{e:ED:bound:gen} \\
		  &p_{u,0} = \overline p_{u,0}, \ \forall u \in \U, \label{e:ED:initial_gen}\\
		  & p_{u,t} - p_{u,t-1} \le \overline R_u u_{u,t-1} + \overline p_{u} v^+_{u,t}, \ \forall   u \in \U, \ t \in [T], \label{e:ED:rup}\\
		& p_{u,t-1} - p_{u,t} \le \underline R_u u_{u,t-1} + \underline p_{u} v^-_{u,t}, \ \forall   u \in \U, \ t \in [T], \label{e:ED:rdown}\\
		& f_{e,t} = -b_e (\theta_{e_h,t} - \theta_{e_t,t}),  \ \forall e \in \E, \ t \in [T], \label{e:ED:DCOPF:PF}\\
		& -\overline f_e \le f_{e,t} \le \overline f_e, \ \forall e \in \E,  \ t \in [T], \label{e:ED:DCOPF:thermalLimit}\\
		& \underline \theta_{i} \le \theta_{i,t}\le \overline \theta_{i}, \ \forall i \in \N, \ t \in [T], \label{e:ED:DCOPF:angleBound}\\
		& -\Delta_{e}  \le \theta_{e_h,t} - \theta_{e_t,t} \le \Delta_{e}  \ \forall e \in \E, \ t \in [T] \label{e:ED:DCOPF:angleDiff}. 
		\end{align}
		\label{prob:e}
\end{subequations}
\vspace{-0.5cm}
\caption{The Unit Commitment and Economic Dispatch Models.}
\label{fig:UCED}
\end{figure}
	
The objective function of the upper level problem (Equations
\eqref{e:obj} - \eqref{e:UC:logic:nonnega}) includes the no-load
costs, the start-up costs, and the costs of the selected supply bids
of each electrical power generating units. Equation \eqref{e:UC:su}
computes the start-up cost $r_{u,t}$ of a generator $u$ for time
period $t$ based on how long $u$ has been offline.  The
expression $u_{u,t} - \sum^h_{n=1} u_{u,t-n}$ is one when generator
$u$ becomes online after it has been turned off for $h$ time
periods. Equation \eqref{e:UC:su:nonnegative} states the nonnegativity
requirement on $r_{u,t}$. Equation \eqref{e:UC:initial_gen_status}
specifies the initial on-off status of each generator. The minimum-up
and -down constraints are specified in Equations \eqref{e:UC:min_up}
and \eqref{e:UC:min_down} respectively. The relationship between the
variables for the on-off, start-up, and shut-down statuses of each
generator is stated in Equation \eqref{e:UC:logic:on_off}. The binary
requirements for logical variables $v^+_{u,t}, v^-_{u,t},u_{u,t}$ are
specified in Equation \eqref{e:UC:logic:nonnega}.

Based on the commitment decisions, the lower-level problem (i.e.,
Equations \eqref{e:ED:obj} - \eqref{e:ED:DCOPF:angleDiff}) decides the
hourly operating schedule of each committed generators in order to
minimize the system production costs. Equation \eqref{e:ED:bal} states
the flow conservation constraints for real power at each bus, using
$e_h$ and $e_t$ to represent the head and tail of $e \in \E$. Equation
\eqref{e:ED:bid} states that the total real power generation of a
generator $u$ is equal to the production of its selected bids.
Equation \eqref{e:ED:bound:bids} constrains the power generation
$s_{b,t}^e$ from bid $b \in \B_u$ to be no more than the submitted amount
$\bar s_b$. Equation \eqref{e:ED:bound:gen} enforces the bound on the
real power generation of each generator. Equation
\eqref{e:ED:initial_gen} specifies the initial generation amount of
each generator, and Equations \eqref{e:ED:rup} and \eqref{e:ED:rdown}
state the ramp-up and -down constraints of each generator. Equation
\eqref{e:ED:DCOPF:PF} captures the DC approximation of the power flow
equations and Equation \eqref{e:ED:DCOPF:thermalLimit} specifies the
thermal limit on each line. Equations \eqref{e:ED:DCOPF:angleBound}
and \eqref{e:ED:DCOPF:angleDiff} state the voltage angle bounds on
each bus and the bounds on the angle difference of two adjacent buses
respectively.


\subsection{The Natural Gas Transmission System}
\label{sec:g}

\begin{table}[t!]\fontsize{9}{9}\selectfont
	\centering
	\caption{Parameters of the gas system} \label{table:param:g}
	{\begin{tabular}{p{0.095\textwidth} p{0.355\textwidth}}
			\toprule
			$\G^g = (\V, \A)$  & Directed graph representing a natural gas transmission network, where $\V$ is a set of junctions, indexed with $j = 1, \cdots, V$, and $\A \subseteq \V \times \V$ is a set of connections, indexed with $a = 1, \cdots, A$ \\%
			$\quad\A_c \subseteq \A$ & Set of compressors\\
			$\quad\A_v \subseteq \A$ & Set of control valves\\
			$\kappa_j$ & Cost of demand shedding at $j \in \V$ \\%
			$ (d^g_{j,t})_{j \in \V}$ & Gas demand profile during $t \in \T$ \\%
			$\underline s^g_{j}, \overline s^g_{j}$ & Lower/Upper limit on natural gas supply at $j \in \V$\\%
			$c_j(\cdot)$ & Cost function for gas supply at $j \in \V$\\%
			$W_a$ & Pipeline resistance (Weymouth) factor of $a \in \A$\\
			$\underline \pi_j, \overline \pi_j$ & Minimum/maximum squared pressure at $j \in \V$\\
			$\underl \alpha^c_a, \overl \alpha^c_a$ & Lower/upper compression ratio of $a \in \A_c$\\
			$\underl \alpha^v_a, \overl \alpha^v_a$ & Lower/upper control ratio of $a \in \A_v$\\
			\bottomrule
	\end{tabular}}{}
\end{table}

\begin{table}[t!]\fontsize{9}{9}\selectfont
	\centering
	\caption{Variables of the gas system }\label{table:var:g}
	{\begin{tabular}{p{0.03\textwidth} p{0.42\textwidth}}
			\toprule
			$s^g_{k,t}$ & Amount of gas supplied by $k \in \K$ during $t \in \T$ \\
			$\pi_{j,t}$ & Pressure squared at $j \in \V$ during $t \in \T$\\
			$\phi_{a,t}$ & Gas flow on $a \in \A$ during $t \in \T$\\
			$l_{j,t}$ & Satisfied gas demand at $j \in \V$ during $t \in \T$\\
			$q_{j,t}$ & Shedded gas demand at $j \in \V$ during $t \in \T$\\
			$\gamma_{j,t}$ & Total amount of gas consumed by the GFPP located at $j \in \N \cap \V$ during $t \in \T$\\
			\bottomrule
	\end{tabular}}{}\end{table}

\begin{figure}[!t]
\begin{subequations}\fontsize{9}{9}\selectfont
		\begin{align}
		\min \ & \sum_{t \in [T]} \sum_{j \in \V}  (\sum_{s \in \mathcal S_j}  c_{j,s} s^g_{s,t} + \kappa_j q_{j,t}) \label{eq:gasobj} \\
		\mbox{s.t.}\  & s^g_{j,t} - l_{j,t} - \gamma_{j,t} = \sum_{a \in \A: a_t = j} \phi_{a,t} - \sum_{a \in \A: a_h = j} \phi_{a,t}, \nonumber\\
		& \hspace{0.25\textwidth}\forall j \in \V, t \in [T], \label{g:bal}\\
		& l_{j,t} = d^g_{j,t} - q_{j,t},  \forall j \in \V, t \in [T],\label{g:load_shed}\\
		& 0 \le q_{j,t} \le d^g_{j,t},  \forall j \in \V, t \in [T],\label{g:load_shed_bound}\\
		& \phi_{a,t} \ge 0,  \forall a \in \A, t \in [T],\label{g:flux_bound}\\
		&  \underline s^g_{j} \le s^g_{j,t} \le \overline s^g_{j}, \forall j \in \V,\ t \in [T],\label{g:bound:supply}\\
		& \underl \alpha^c_a \pi_{a_h,t} \le \pi_{a_t,t}  \le  \overl \alpha^c_a \pi_{a_h,t}, \forall a  \in \A_c,\ t \in [T], \label{g:compressors}\\
		& \underl \alpha^v_a \pi_{a_h, t}  \le \pi_{a_t,t}  \le  \overl \alpha^v_a \pi_{a_h,t}, \forall a \in \A_v,\ t \in [T], \label{g:valves}\\
		&  \pi_{a_h, t}  - \pi_{a_t,t}  = W_a \phi_{a,t}^2,\forall a \in \A \setminus (\A_v \cup \A_c),\ t \in [T], \label{g:Weymouth}\\
		&\underline \pi_j  \le \pi_{j, t}  \le \overline \pi_{j},\forall j \in \V,\ t \in [T] \label{g:bound:pi} \\
                & s^g_{j,t} = \sum_{s \in \mathcal S_j} s^g_{s,t} \label{g:link}
		\end{align}
		\label{prob:g}
\end{subequations}
\vspace{-0.5cm}
\caption{The Natural Gas Transmission Model.}
\label{fig:gas}
\end{figure}

Tables \ref{table:param:g} and \ref{table:var:g} specify the
parameters and variables of the steady-state natural gas model, which
is given in Figure \ref{fig:gas}. The modeling is similar to those in
\cite{bent2018joint,sanchez2016convex,borraz2016convex} and uses the
Weymouth equation to capture the relationship between pressures and
flux.  The flux conservation constraint is given in Equation
\eqref{g:bal}, where $a_h$ and $a_t$ represent the head and tail of $a
\in \A$. Equation \eqref{g:load_shed} determines the demand served at
each junction: It captures the amount of gas load shedding which must
be nonnegative and cannot exceed the demand at the corresponding
junction (Equation \eqref{g:load_shed_bound}). The model assumes that
gas flow directions are predetermined and Equation
\eqref{g:flux_bound} enforces the sign of gas flow variables, i.e., it
constrains $\phi_{a,t}$ to be nonnegative. Equation
\eqref{g:bound:supply} specifies the upper and lower limits of natural
gas supplies. The change in pressure through compressors and control
valves are formulated in Equations \eqref{g:compressors} and
\eqref{g:valves} and the model use a single compressor machine
approximation as in prior work.  The steady-state physics of gas flows
is formulated with the Weymouth equation in Equation
\eqref{g:Weymouth}. Equation \eqref{g:bound:pi} states the bounds on
nodal pressures. Equation \eqref{g:Weymouth} can be convexified using
the second-order cone relaxation from \cite{borraz2016convex}:
$
\pi_{a_h, t} - \pi_{a_t,t} \ge W_a \phi_{a,t}^2.
$
This relaxation is very tight \cite{borraz2016convex}.

\label{sec:g:cost}


\label{rem:price}
When the gas system is not congested, the price of natural gas is
relatively stable. However, during congestion and when some loads are
being shedded, natural gas prices increase sharply. The cost of gas in
the objective function captures this behavior: For a junction $j$, it
is specified with an almost-linear piecewise linear function for
production and a high penalty cost $\kappa_j$ for gas shedding. To be
specific, let $\mathcal S_j$ be a set of non-overlapping intervals
covering $[0, \overline s^g_j]$, each with a distinct slope $c_{j,s}$
satisfying $c_{j,s} \le c_{j,s+1}$ whenever $s,s+1 \in S_j$.  Define
an auxiliary nonnegative variable $s^g_{s,t}$ that represents the
amount of gas supply from $s \in S_j$ at time $t$. The objective
function is then stated as
\begin{equation}
\sum_{t \in [T]} \sum_{j \in \V}  (\sum_{s \in \mathcal S_j}  c_{j,s} s^g_{s,t} + \kappa_j q_{j,t}).    
\nonumber
\end{equation}
The model also includes constraint \eqref{g:link} to link the gas
variable at junction $j$ with the auxiliary variables.

\subsection{Physical and Economic Couplings}
\label{sec:c}

GFPPs are the physical and economic interface between the electrical
power and gas networks. This section first describes the resulting
coupling constraints before describing how the natural gas zonal
prices are computed. Tables \ref{table:param:c} and \ref{table:var:c}
describe the parameters for the coupling.

\begin{table}[t]\fontsize{9}{9}\selectfont
	\centering
	\caption{Parameters for the Electricity and Gas Coupling.} \label{table:param:c}
	{\begin{tabular}{p{0.1\textwidth} p{0.35\textwidth}}
			\toprule
			$\{H_{u,i}\}_{i = 0,1,2}$ & Coefficients of the heat rate curve of $u \in \U^g$\\%
			$\alpha_u$ & Maximum allowable percentage of the expense on natural gas over its marginal bid price for $u \in \U^g$\\
			$\mathcal K $ & Set of pricing zones, indexed with $k = 1, \cdots, K$\\
			$\quad \V(k) $ & Set of junctions that belong to $k \in \K$\\
			\bottomrule
	\end{tabular}}{}
\end{table}

\begin{table}[t]\fontsize{9}{9}\selectfont
	\centering
	\caption{ Variables for the Electricity and Gas Coupling.}\label{table:var:c}
	{\begin{tabular}{p{0.03\textwidth} p{0.42\textwidth}}
			\toprule
			$w_{b,t}$ &  1 if $b \in \B_u$ of $u \in \U$ is selected during $t \in \T$, 0 otherwise \\
			$\rho_{u,t}$ & Price of marginally selected bid of $u \in \U^g$ during $t \in \T$\\
			$\psi_{k,t}$ & Zonal price of natural gas in $k \in \K$ during $t \in \T$\\
			\bottomrule
	\end{tabular}}{}\end{table}

The physical couplings between $\G^e$ and $\G^g$ can be formulated as follows ($t \in [T], j \in \N \cap \V$):
\begin{align}
\gamma_{j,t} = \sum_{u \in \U(i) \cap \U^g} H_{u,2} p_{u,t}^2 + H_{u,1} p_{u,t} + H_{u,0}.	\label{c:physical}
\end{align}
The real power generation $\boldsymbol{p}$ of a GFPP induces a demand
$\boldsymbol{\gamma}$ in the natural gas system. Equation
\eqref{c:physical} specifies the relationship between the real power
generation of a GFPP and the amount of natural gas needed for the
generation. In the equation, this relationship is approximated by a
quadratic heat-rate curve, whose coefficients are given as $H_u$. The
equation can be convexified like the Weymouth equation.

Since the level of power generation of the GFPPs determines the load
in the gas system, the physical coupling also affects the natural gas
prices. The price formation of natural gas, in turn, governs the
profitability of GFPPs, which submit bids before the realization
of gas prices. To capture these economic realities, the model
introduces binary variables of the form $w_{b,t} \in \{0,1\}$ for each
bid $b$ of a GFPP: Variable $w_{b,t}$ indicates whether bid $b$ is
selected during time period $t$. Equation \eqref{e:ED:bid} is then
replaced by the following constraints (for all $t \in [T]$):
\begin{subequations}
    \begin{align}
      &\rho_{u,t} = \sum_{b \in [B_u-1]} \beta_b (w_{b,t} - w_{b+1,t}) + \beta_{B_u} w_{B_u,t},  \forall u \in \U^g \label{c:bid_price},\\
&  0 \le s^e_{b,t} \le \overline {s}_{b},  \forall b \in \B_u, \ u \in \U \setminus \U^g, \label{e:UC:bound:bids:nongfpp} \\
			&  0 \le s^e_{b,t} \le \overline {s}_{b}w_{b,t}, \forall b \in \B_u, \ u \in \U^g \label{e:UC:bound:bids:gfpp}\\
&  w_{b,t} \le u_{u,t}, \forall b \in \B_u, \ u \in \U^g,  \label{e:UC:logic:bid1}\\
			&  \overline s_b w_{b+1,t} \le s_{b,t},  \forall b \in [1, B_u-1]_{\mathbb Z}, \ u \in \U^g \label{e:UC:logic:bid2}.
    \end{align}
    \label{e:UC:bidsBound2}
\end{subequations}

\noindent
Equations \eqref{e:UC:bound:bids:nongfpp} and
\eqref{e:UC:bound:bids:gfpp} are bound constraints for the bids
submitted by the non-GFPPs and GFPPs
respectively. Equation \eqref{e:UC:bound:bids:gfpp} ensures that the
indicator variable $w_{b,t}$ is one whenever bid $b$ is used for time
period $t$ (i.e., $s^e_{b,t} > 0$). Equation \eqref{e:UC:logic:bid1}
states that the bid of a generator can be selected only when it is
committed and Equation \eqref{e:UC:logic:bid2} ensures that the
$(b+1)^{\mbox{th}}$ bid is selected only if the bid $b$ is fully
used. Accordingly, Equation \eqref{c:bid_price} states that
$\rho_{u,t}$ is the maximum/marginal bid price of GFPP $u \in \U^g$
among its currently selected bids.

The economic coupling between the electricity and gas networks is
enforced by {\em bid-validity constraints} that ensure that the
marginal costs of producing electricity by GFPPs are lower than their
marginal bid prices. Although the natural gas system is operated in a
decentralized manner, the zonal price of natural gas $\boldsymbol\psi$
can be modeled as a function $g$ of the market supply and demand,
i.e., as a function of the binary and continuous variables of Problems
\eqref{prob:e} and \eqref{prob:g}, which are denoted by $\z$ and $\x$.
Under this assumption, the bid validity constraints can be expressed
as follows (for all $t \in [T]$):
\begin{subequations}
  \begin{align}
                & \psi = g(z, x), \label{c:price}\\
		& \alpha_u \rho_{u,t}+M (1-u_{u,t}) \ge \left[2p_{u,t}H_{u,2} + H_{u,1} \right]\psi_{k,t},\nonumber\\
		& \hspace{0.13\textwidth}\forall k \in \K,  i \in  \V(k),  u \in \U(i) \cap \U^g.  \label{c:bid_val}
		\end{align}
		\label{c:economic}
\end{subequations}
\noindent
They capture the fact that, when the realized natural gas price
\[
\left[2p_{u,t}H_{u,2} +  H_{u,1} \right] \psi_{k,t}
\]
for generating one additional unit of real power by GFPP $u$ is
greater than its marginal bid price $\rho_{u,t}$, GFPP $u$ is not
profitable.  This situation arises because GFPP $u$ submits its bids
before the realization of $\psibold$. The bid validity constraint is
expressed in Equation \eqref{c:bid_val} and ensures that only
profitable GFPPs are committed. The bid validity constraints use the
realized zonal gas prices from Equation \eqref{c:price} and $M$
denotes a big-M value set to the maximum natural gas price (e.g.,
\$200 per mmBtu) multiplied by $\left[2 \overline p_{u} H_{u,2} +
  H_{u,1} \right]$.

It remains to specify how to compute the zonal gas prices, i.e., the
function $g$ in Equation \eqref{c:price}. The UCGNA assumes that the
nodal natural gas price at each junction $j$ is given by the marginal
cost of supplying natural gas at $j$. This marginal cost is the dual
solution associated with the corresponding flux conservation
constraint in Problem \eqref{prob:g}. The zonal natural gas prices
$\psibold$ are then computed by averaging the nodal natural gas prices
of a subset of junctions in the zone. Therefore, the zonal natural gas
price $\psibold$ are given by linear functions of the dual solution to
Problem \eqref{prob:g}.

Note that, by construction, the natural gas zonal prices $\psibold$
under normal operating conditions are given by the almost linear part
of objective \eqref{eq:gasobj}. However, when the gas network is
congested and load needs to be shed, the zonal prices increase
sharply due to the high penalty cost $\kappa_j$. As a result, the
resulting model closely captures the behavior of the market during the
2014 polar vortex. Note also that the model does not shed the demand
of the GFPPs. The model assumes that GFPPs buy natural gas at any cost
to meet its commitment obligation. Once again, this captures the 2014
Polar Vortex situation where GFPPs were encouraged to buy the natural
gas from the spot market at any cost for the sake of the power system
reliability \cite{FERC2015}.


\section{Reformulation of the UCGNA}
\label{sec:approxi}

This section shows how the UCGNA can be expressed as a MISOCP.  Let
variable subscripts $p$ and $g$ respectively denote the power and the
gas systems. Let $\z_p$ and $\x_p$ respectively denote the vector of
binary and continuous variables of the power system (i.e., Problem
\eqref{prob:e}) and let $\x_g$ be the vector of continuous variables
of the gas system (i.e., Problem \eqref{prob:g}).  The UCGNA can be
stated as a trilevel program:
	
\begin{subequations}\fontsize{9}{9}\selectfont
			\begin{align}
			\min_{\underset{\z_p \in \{0,1\}^m}{\x_p \ge 0,  \y_g }} \quad & c_p^T \x_p + h^T \z_p \label{TL:obj:leader}\\
			\mbox{s.t. } & \;\;\; \z_p \in \mathcal Z, \label{TL:1level}\\
			& \begin{array}{lcl}
				(\x_p, \y_g) = &  \underset{\x_p \ge 0, \y_g }{\mbox{argmin}} & c_p^T \x_p\\ 
				& \mbox{s.t.} &  A \x_p + B \z_p \ge b, \\
				& &	\y_g \in \mbox{Dual sol. of } \eqref{TL:g},\\
				\end{array} \label{TL:lower}\\
			& \;\;\; E \y_g + M \z_p \ge h \label{TL:additional}
			\end{align}
			\label{prob:TL}
\end{subequations}
\noindent
where $\mathcal Z$ denotes the feasible region of the unit commitment
problem (i.e., Equations \eqref{e:UC:su}-\eqref{e:UC:logic:nonnega}),
the third level problem is defined as
\begin{equation}
  \underset{\x_g \in \K}{\min}  \ c_g^T \x_g : D_p \x_p + D_g \x_g \ge d, \label{TL:g}
\end{equation}
and $\K$ is the proper cone denoting the domain of $\x_g$.

The first-level problem (i.e., Equations \eqref{TL:obj:leader} and
\eqref{TL:1level}) formulates the unit-commitment problem (i.e.,
Equations \eqref{e:obj}-\eqref{e:UC:logic:nonnega} and Equation \eqref{e:UC:bidsBound2}). The
unit-commitment decisions $\z_p$ from the first-level problem are then
plugged into the second-level problem, which formulates the economic
dispatch problem (i.e., Equations
\eqref{e:ED:obj}-\eqref{e:ED:DCOPF:angleDiff}) and decides the hourly
operating schedule of committed generating units. Then, the
third-level problem (i.e., Problem \eqref{TL:g}) formulates the
natural gas problem (i.e., Problem \eqref{prob:g} and Equation
\eqref{c:physical}) and determines the resulting nodal prices for
natural gas based on the dual solution $\y_g$ of the economic dispatch
decisions.

Equations \eqref{TL:obj:leader}, \eqref{TL:1level}, and
\eqref{TL:lower} capture the current operating practice of the power
system. The first level captures the commitment decisions that are
taken first without consideration of the gas network. The second and
third levels implement a Stackelberg game, where the dispatch
decisions of the electricity system are followed by those of the
natural gas network. {\em The novelty in the UCNGA is the bid-validity
  constraint \eqref{TL:additional}, which corresponds to Equation
  \eqref{c:bid_val}}: It ensures that only profitable GFPPs are
selected in the first level and uses the dual variables of the
third-level problem to do so, allowing the unit-commitment problem to
anticipate the zonal prices of natural gas.

The following theorem, whose proof is in Appendix \ref{appendix:pf},
shows that the tri-level problem can be reformulated as a single-level
mathematical program. The proof uses strong duality on the third-level
problem and a lexicographic optimization to merge the second and third
levels.
	
\begin{theorem}
  Problem \eqref{prob:TL} can be asymptotically approximated by the following
  mathematical program:

  \begin{subequations}\fontsize{9}{9}\selectfont
\begin{align}
\min \ &  \alpha h^T \z_p + \alpha c_p^T \x_p  + (1-\alpha) c_g^T  \x_g  \label{our:obj}\\
\mbox{s.t.} \ &	\z_p \in \mathcal Z, \\
& A \x_p + B \z_p \ge b,\label{our:leader}\\
& D_p \x_p + D_g \x_g \ge d,\\
& \y_p^T (b - B \z_p) +  \y_g^T d  \ge \alpha c_p^T \x_p  +(1-\alpha)  c_g^T \x_g, \label{our:strong}\\
& \y_g^T D_g  \preceq_{\K^*} (1-\alpha)  c_g^T,\\
& \y_p^T A + \y_g^T D_p  \le \alpha c^T_p,\\
& 	\frac{1}{1-\alpha} E  \y_g  + M \z_p \ge h,\\
&\x_p \ge 0, \x_g \in \K, \y_p \ge 0, \y_g \ge 0, \\
&\z_p \in \{0,1\}^m,
\end{align}
\label{prob:our}
\end{subequations}
for some $\alpha \in (0,1)$. Moreover, when $\alpha \rightarrow 1$, 
the optimal solution of Problem \eqref{prob:our} converges to
the optimal solution of Problem \eqref{prob:TL}.
\label{theo:TL}
\end{theorem}

\noindent
Observe that Problem \eqref{prob:our} has a bilinear term of $\y_p^T B
\z_p$ in Equation \eqref{our:strong}. Assuming that $\y$ has an upper
bound of $\overline \y$, this term can be rewritten using an exact
McCormick relaxation to produce a MISOCP.

\begin{remark}\label{rem:2}
Problem \eqref{prob:our} is best viewed as a ``standard'' MISOCP to
which a constraint on the dual variables of its inner-continuous
problem has been added.  The ``standard'' MISOCP optimizes the joint
electricity and natural gas problem

\begin{subequations}\fontsize{9}{9}\selectfont
\begin{align}
\min \ &   \alpha h^T \z_p+ \alpha c_p^T \x_p  +  (1-\alpha) c_g^T \x_g   \\
\mbox{s.t.} \ & 	\z_p \in \mathcal Z, \\
	& A \x_p + B \z_p \ge b,\\
	& D_p \x_p + D_g \x_g \ge d,\\
	& \x_p \ge 0, \x_g \in \K, \z_p \in \{0,1\}^m,
\end{align}
\label{prob:p}
\end{subequations}
\noindent
and the additional constraints
\[
\frac{1}{1-\alpha} E y_g + M z_p \ge h.
\]
on the dual variables $(y_p, y_g)$ of its inner continuous problem
capture the bid validity.
\end{remark}

\section{Solution Approach}
\label{sec:sol_approach}

This section briefly sketches how the MISOCP is solved.  Problem
\eqref{prob:our} can be reformulated as

\begin{subequations}\fontsize{9}{9}\selectfont
	\begin{align}
	\underset{\z_p \in \mathbb{B}^n}{\min} \ & \alpha h^T \z_p +  f(\z_p) \\
	\mbox{s.t.} \ & \z_p \in \mathcal Z.
	\end{align}
\end{subequations}
where
\begin{subequations}\fontsize{9}{9}\selectfont
\begin{align}
f(\z_p) = \ & \min  \alpha c_p^T \x_p  + (1-\alpha) c_g^T  \x_g \\
\mbox{s.t.} \ &	A \x_p + B \z_p \ge b, \\
& D_p \x_p + D_g \x_g \ge d,\\
& \y_p^T (b - B \z_p) +  \y_g^T d  \ge \alpha c_p^T \x_p  +(1-\alpha)  c_g^T \x_g, \\
& \y_g^T D_g  \preceq_{\K^*} (1-\alpha)  c_g^T,\\
& \y_p^T A + \y_g^T D_p  \le \alpha c^T_p,\\
& 	\frac{1}{1-\alpha} E  \y_g  + M \z_p \ge h,\\
&\x_p \ge 0, \x_g \in \K, \y_p \ge 0, \y_g \ge 0.
\end{align}
\label{prob:Benders}
\end{subequations}

\noindent
The implementation applies a Benders decomposition on this formulation
to solve Problem \eqref{prob:our}. Moreover, the dual of Problem
\eqref{prob:Benders} has a special structure that can be exploited by
the dedicated Benders decomposition from \cite{ByeonTheory2018}. The
idea is to decompose the dual of Problem \eqref{prob:Benders} into two
more tractable problems. The extreme points and rays of these
subproblems can be used to find the (feasibility and optimality)
Benders cuts of Problem \eqref{prob:Benders}.  The solution method
also uses the acceleration schemes from \cite{fischetti2010note,Fischetti2017RedesigningBD} which
normalize the rays $\hat y$ and perturb $\hat z_p$. The solution
method also obtains feasible solutions periodically (e.g., every 30
iterations) heuristically by turning off violated generators. Finally,
the solution method applies a preprocessing step to eliminate some
invalid bids. It exploits the fact that the natural gas prices without
the GFPP load gives a lower bound on the natural gas zonal
prices. Therefore, the implementation solves Problem \eqref{prob:g}
with no GFPPs, i.e., $\gamma_{j,t} = 0$ for all $j \in \V, t \in
[1,T]_{\mathbb Z}$. Those bids violating the bid-validity constraint
with regard to these zonal prices are not considered further.
			
\section{Description of the Data Sets} \label{sec:data}

The UCGNA model is evaluated on the gas-grid test system from
\cite{bent2018joint}, which is representative of the natural gas and
electric power systems in the Northeastern United States. This test
system is composed of the IEEE 36-bus NPCC electric power system
\cite{allen2008combined} and a multi-company gas transmission network
covering the Pennsylvania-To-Northeast New England area in the United
States \cite{bent2018joint}. The data for the test system can be found
online at \url{https://github.com/lanl-ansi/GasGridModels.jl}.

The test system consists of 91 generators of various types (e.g.,
hydro, gas-fueled, coal-fired, nuclear, etc.).  The unit-commitment
data for these generators (e.g., generator offer curves including
start-up and no-load costs and operational parameters such as minimum
run time) was obtained from the RTO unit commitment test system
\cite{krall2012rto}. Each generator in the gas-grid test system is
assigned the unit commitment data adapted to its fuel-type and megawatt
capacity.

The gas-grid test case consists of two natural gas pricing zones:
Transco Zone 6 non NY and Transco Leidy Line. The Transco Leidy Line
represents the natural gas prices in the Marcellus Shale production
area, which has a wealth of natural gas. On the other hand, the
Transco Zone 6 non NY represents the natural gas prices near
consumption points. Therefore, a large difference in prices between
these two pricing zones implies a scarcity of transmission capacities
between these two points. During normal operations, the average
natural gas prices in the Transco Zone 6 non NY and the Leidy Line are
around \$3/mmBtu and \$1.5/mmBtu respectively. The slopes $c_{j,s}$ at
junction $j \in \V$ (see Section \eqref{sec:g:cost}) are chosen to be
around these numbers. The penality cost for load shedding $\kappa_aj$
is set as \$130/mmBtu for all junctions. The results are given for a
single time-period (i.e., $T = 1$).


\section{Case Study}
\label{sec:case_study}

This section analyzes, under various operating conditions, the
behavior of the UCGNA on the realistic test system described in
Section \ref{sec:data}. The results are compared with current
practices. The case study varies the level of stress on both the
electrical power and gas systems. For the electrical power system, the
load is uniformly increased by 30\% and 60\%. For the gas system, the
load is uniformly increased by 10\% up to 130\%. Parameters $\eta_p$
and $\eta_g$ respectively represent the stress level imposed on the
electrical power and gas systems. In the results, (A) denotes existing
practices and (B) the UCGNA model. Solutions for (B) are obtained with
a wall-clock time limit of 1 hour, while solutions for (A) is obtained
by the following procedure:
\begin{enumerate}
    \item[(i)] Solve the power model (i.e., Problem \eqref{prob:e});
    \item[(ii)] Retrieve the demand of GFPPs using Equation
      \eqref{c:physical} and plug it into the gas model (i.e., Problem
      \eqref{prob:g});
    \item[(iii)] Solve the gas model and compute the natural gas 
      zonal prices using the dual values associated with the flux
      conservation constraints;
    \item[(iv)] Based on the zonal prices, determine the set of GFPPs
      violating the bid-validity constraint (i.e., Equation
      \eqref{c:bid_val}) and compute the loss of such GFPPs by
      multiplying the violation, i.e., the difference between the marginal gas price and the
      marginal bid price, with the scheduled amount of power
      generation.
\end{enumerate}

\begin{figure}[!t]
	\centering
	\subfloat[System costs (A).\label{fig:pStress:1:cost:A}] {\includegraphics[width=0.22\textwidth]{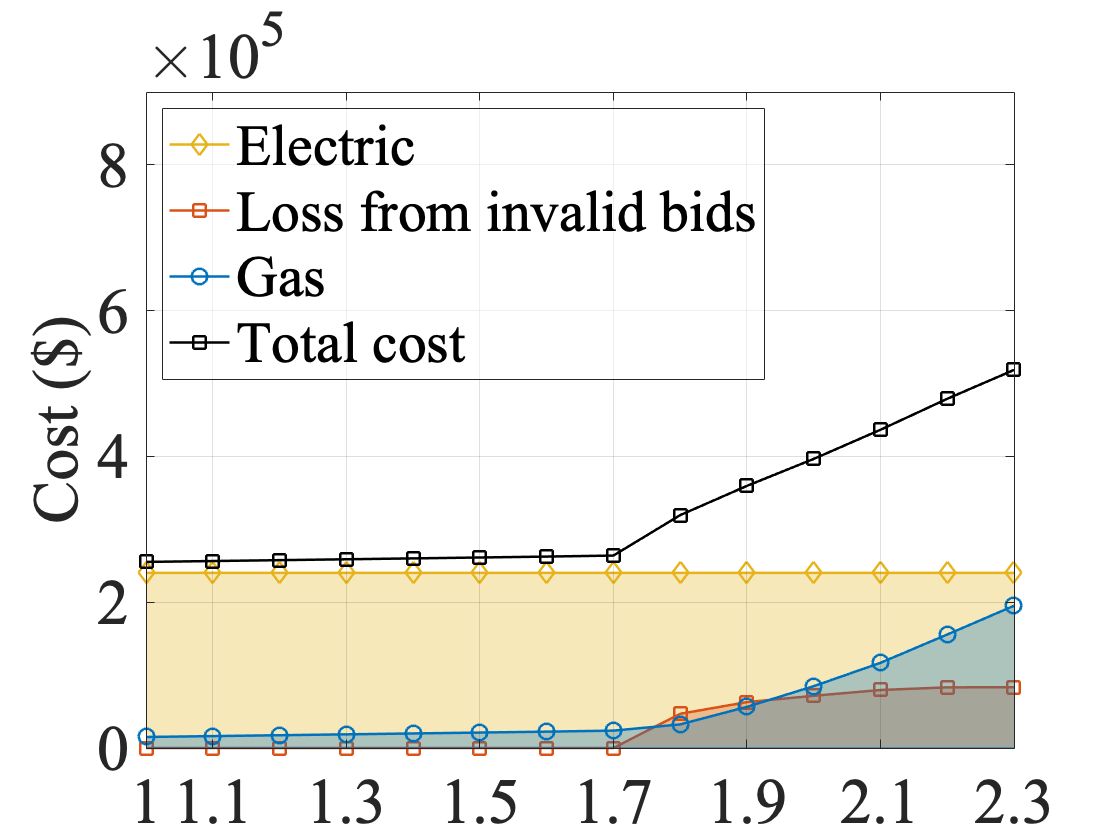}}
	\subfloat[System costs (B).\label{fig:pStress:1:cost:B}] {\includegraphics[width=0.22\textwidth]{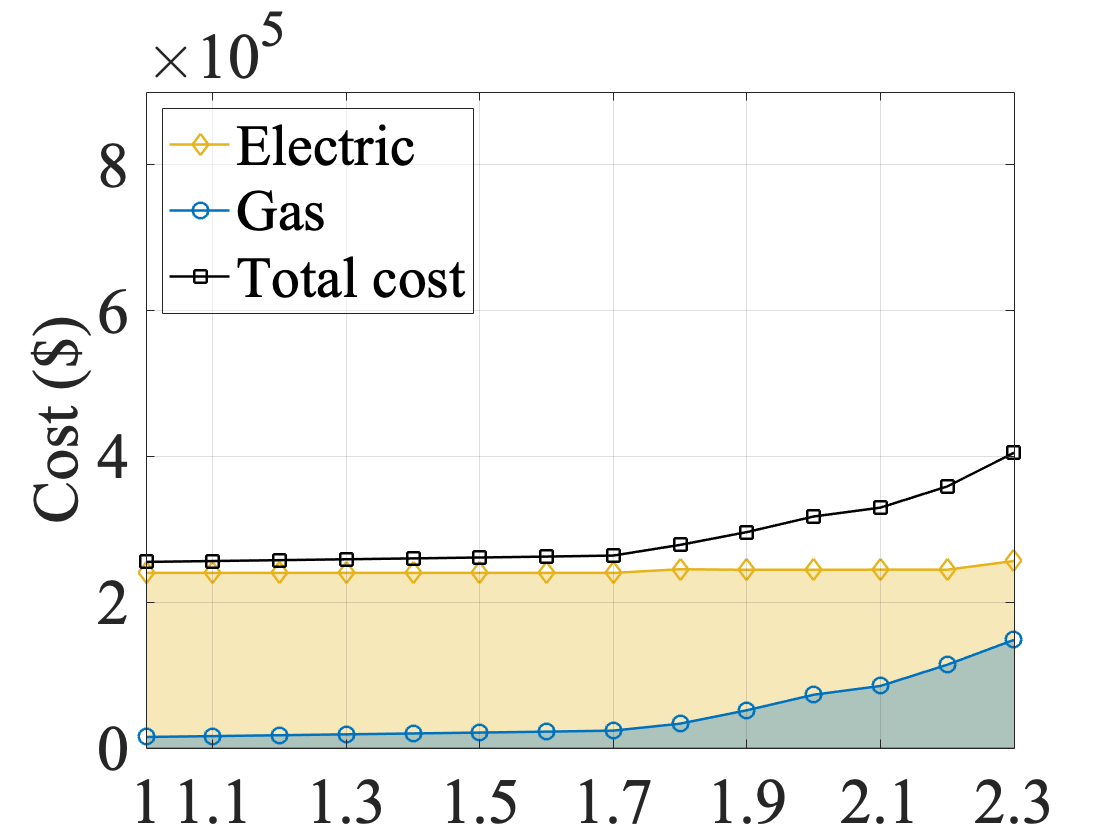}}\\
	\subfloat[Natural gas prices (A).\label{fig:pStress:1:price:A}] {\includegraphics[width=0.22\textwidth]{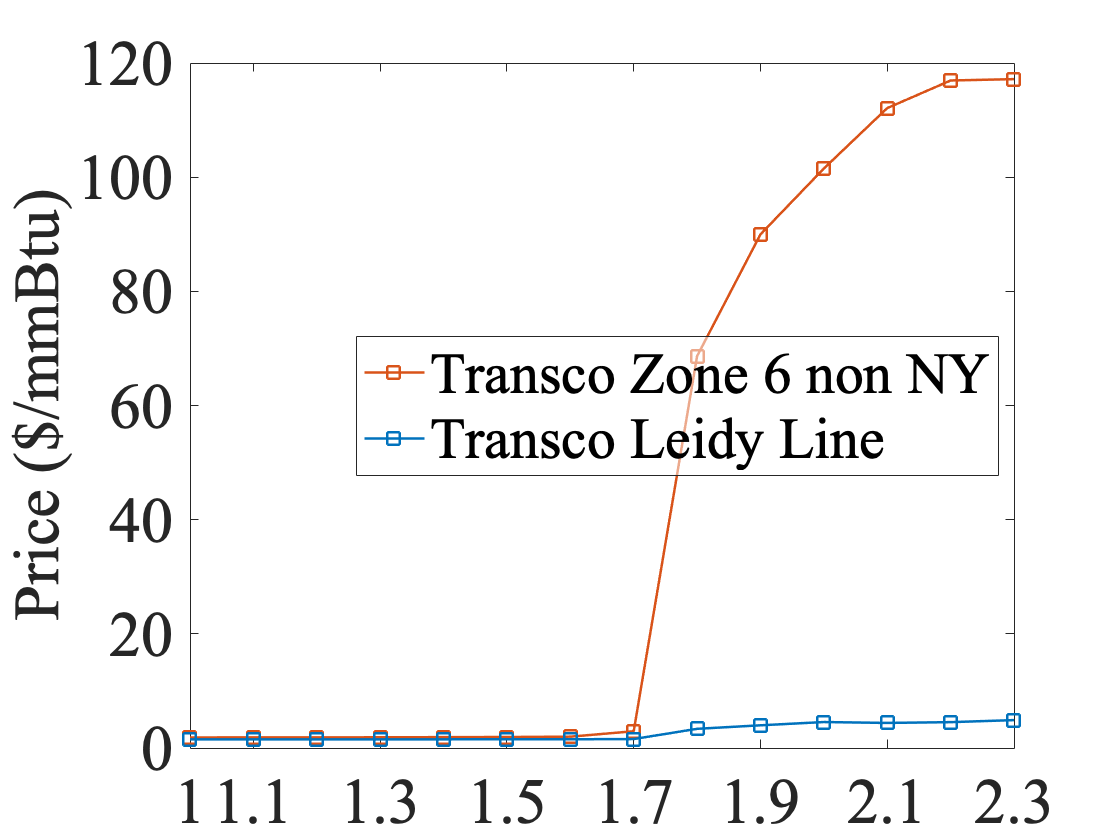}}
	\subfloat[Natural gas prices (B).\label{fig:pStress:1:price:B}] {\includegraphics[width=0.22\textwidth]{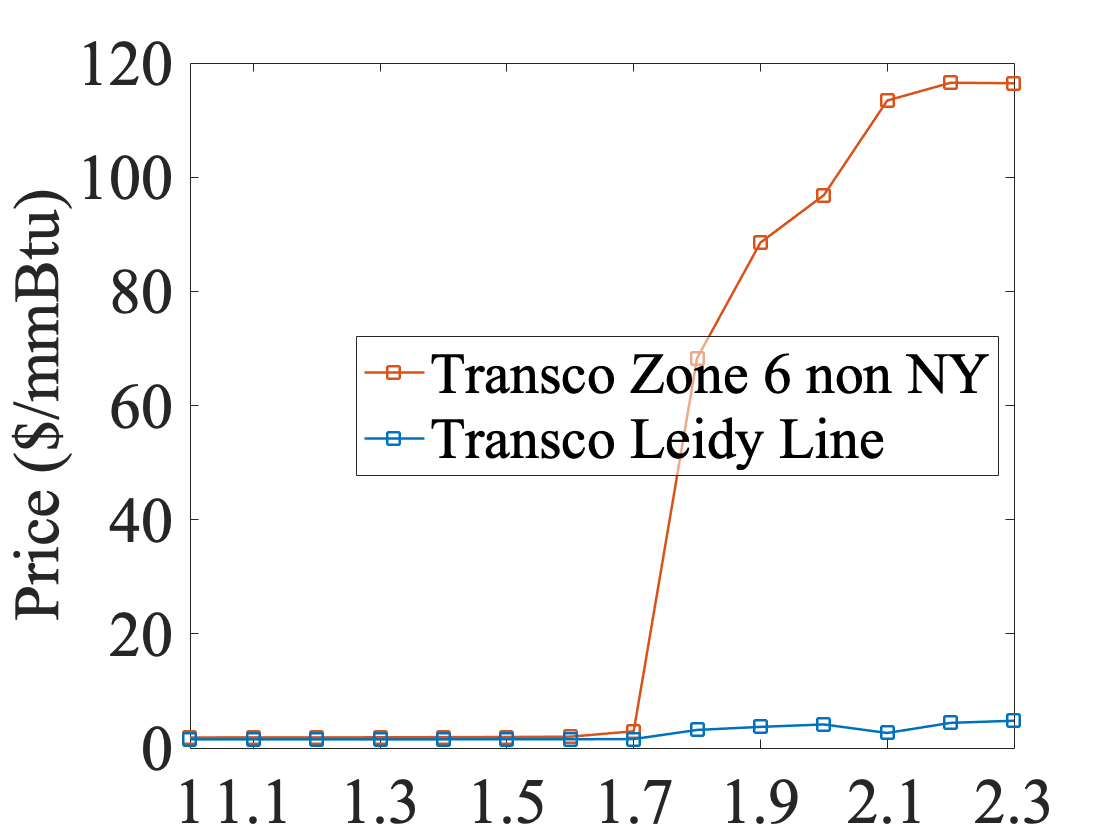}}
	\caption{Results for the Normal Operating Conditions of the Electrical Power System ($\eta_p = 1$).} \label{fig:pStress:1}
\end{figure}

\begin{figure}[!t]
	\centering
	\subfloat[System costs (A).\label{fig:pStress:13:cost:A}] {\includegraphics[width=0.22\textwidth]{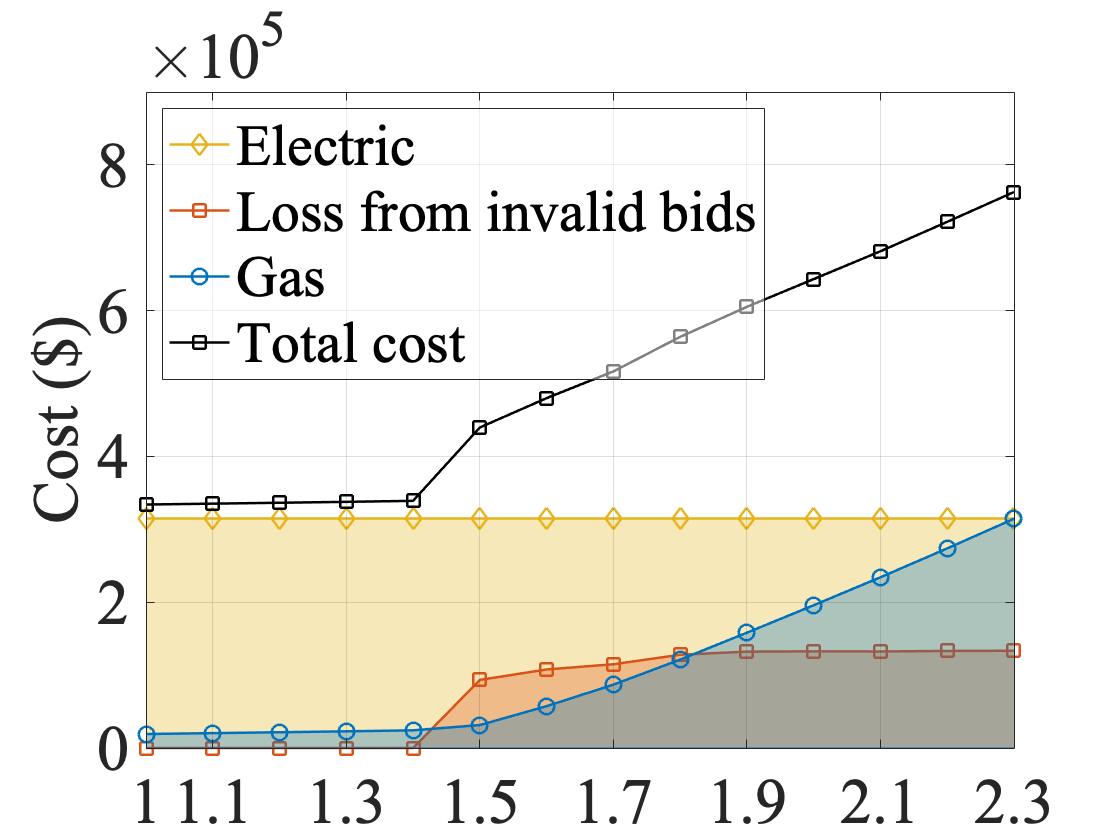}}
	\subfloat[System costs (B).\label{fig:pStress:13:cost:B}] {\includegraphics[width=0.22\textwidth]{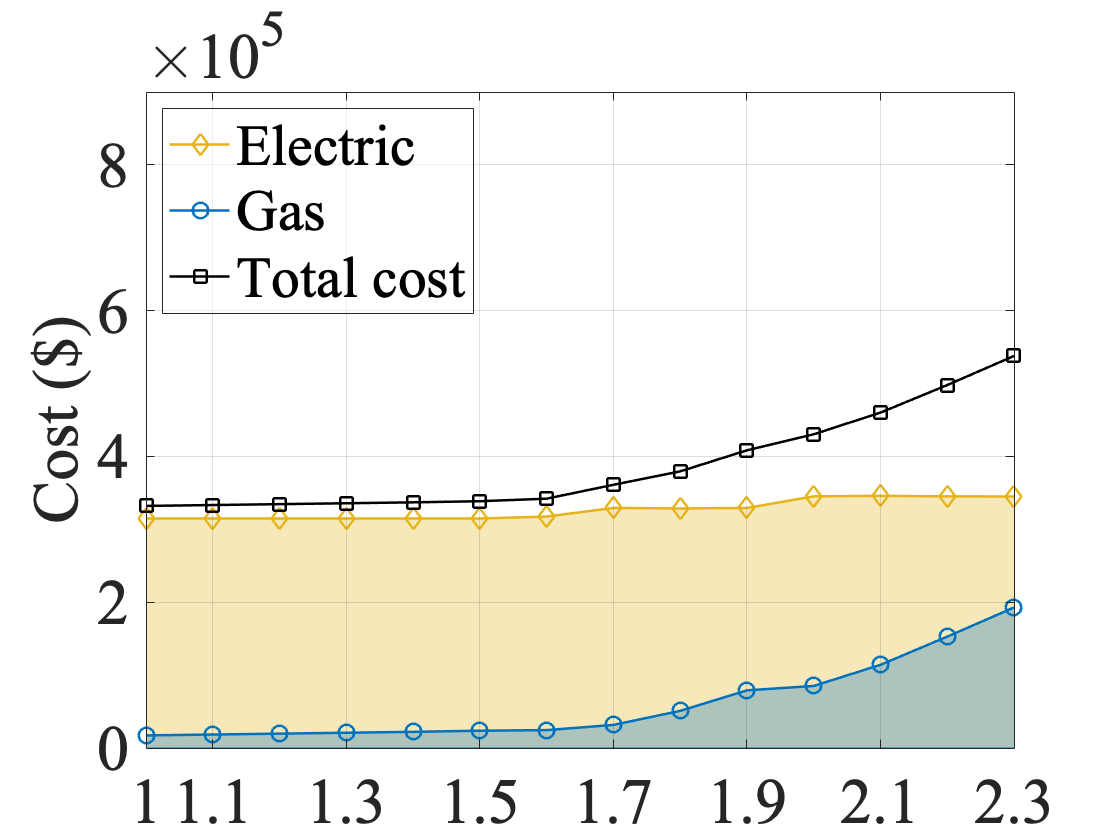}}\\
	\subfloat[Natural gas prices (A).\label{fig:pStress:13:price:A}] {\includegraphics[width=0.22\textwidth]{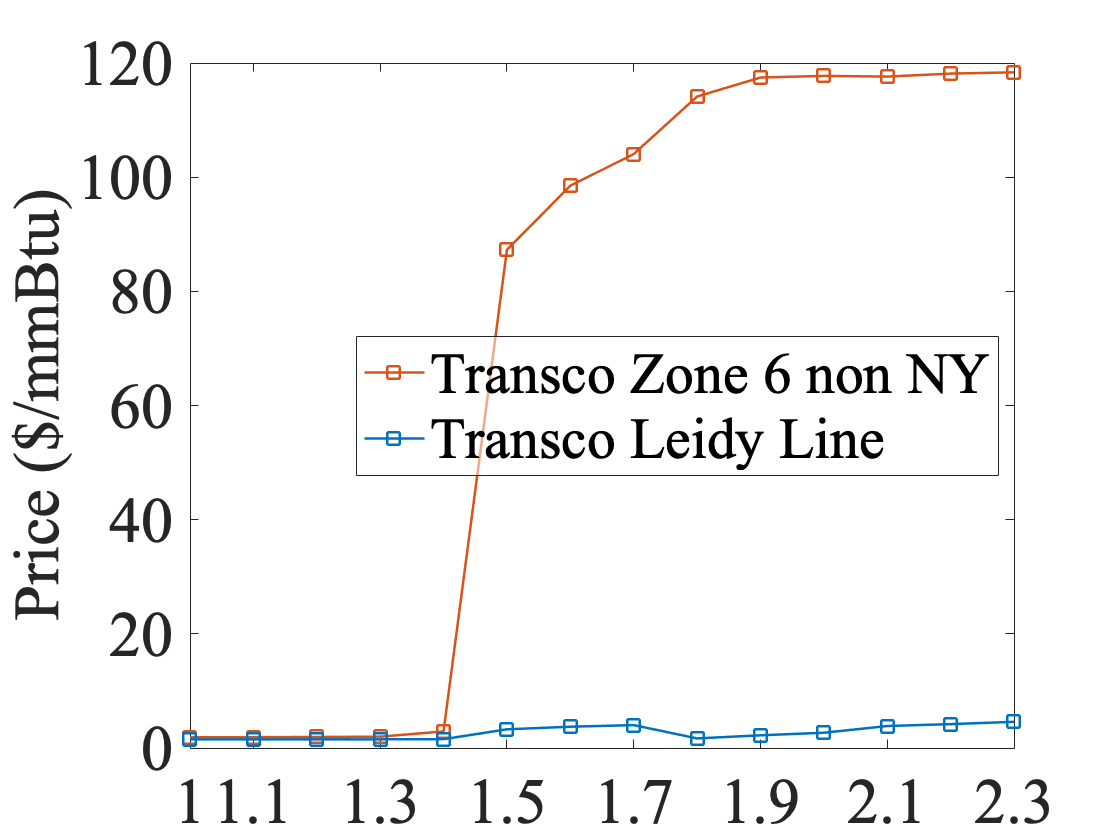}}
	\subfloat[Natural gas prices (B).\label{fig:pStress:13:price:B}] {\includegraphics[width=0.22\textwidth]{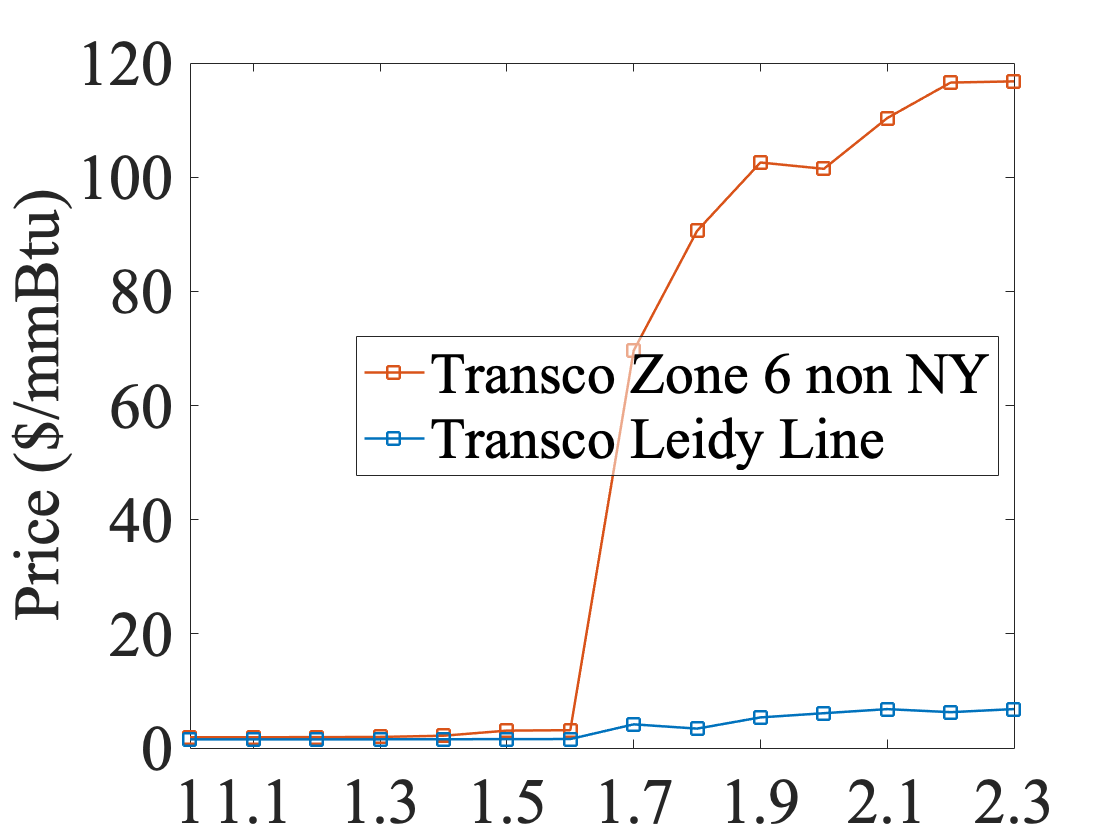}}
	\caption{Results for the Stressed Electrical Power System ($\eta_p = 1.3$).} \label{fig:pStress:13}
	\vspace{-0.25cm}
\end{figure}

\begin{figure}[!t]
	\centering
	\subfloat[System costs (A).\label{fig:pStress:17:cost:A}] {\includegraphics[width=0.22\textwidth]{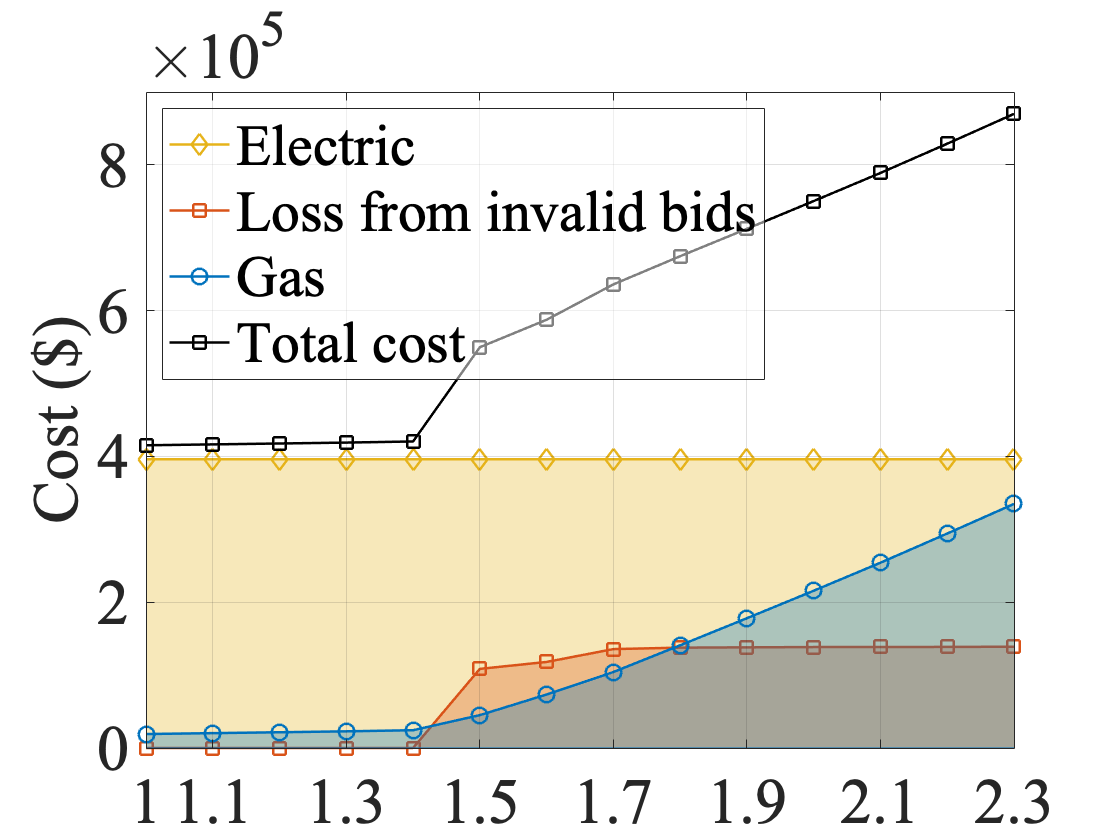}}
	\subfloat[System costs (B).\label{fig:pStress:17:cost:B}] {\includegraphics[width=0.22\textwidth]{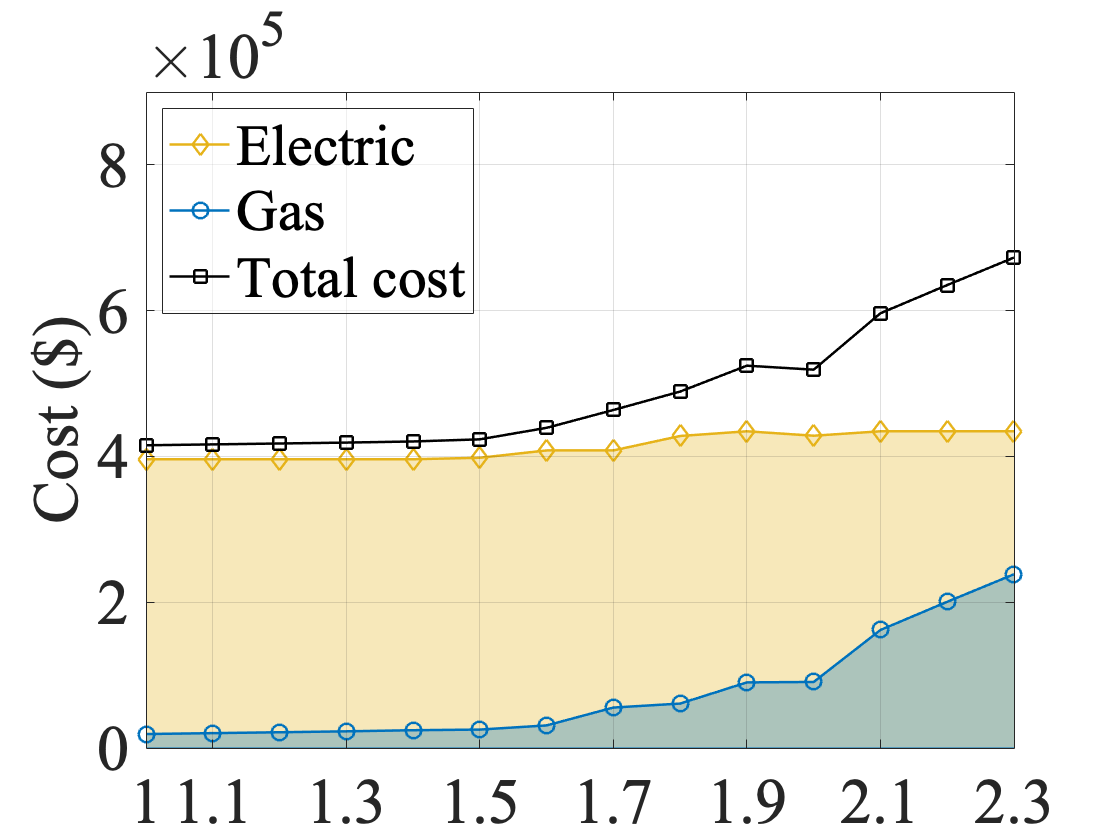}}\\
	\subfloat[Natural gas prices (A).\label{fig:pStress:17:price:A}] {\includegraphics[width=0.22\textwidth]{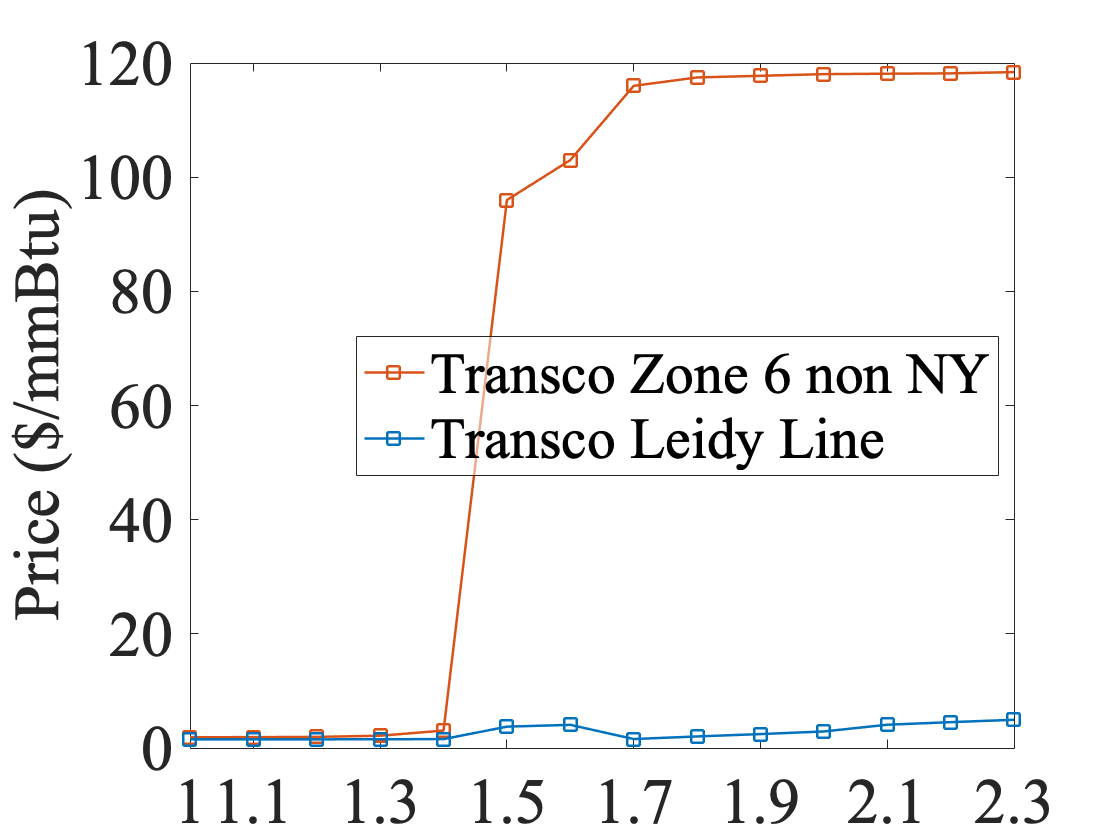}}
	\subfloat[Natural gas prices (B).\label{fig:pStress:17:price:B}] {\includegraphics[width=0.22\textwidth]{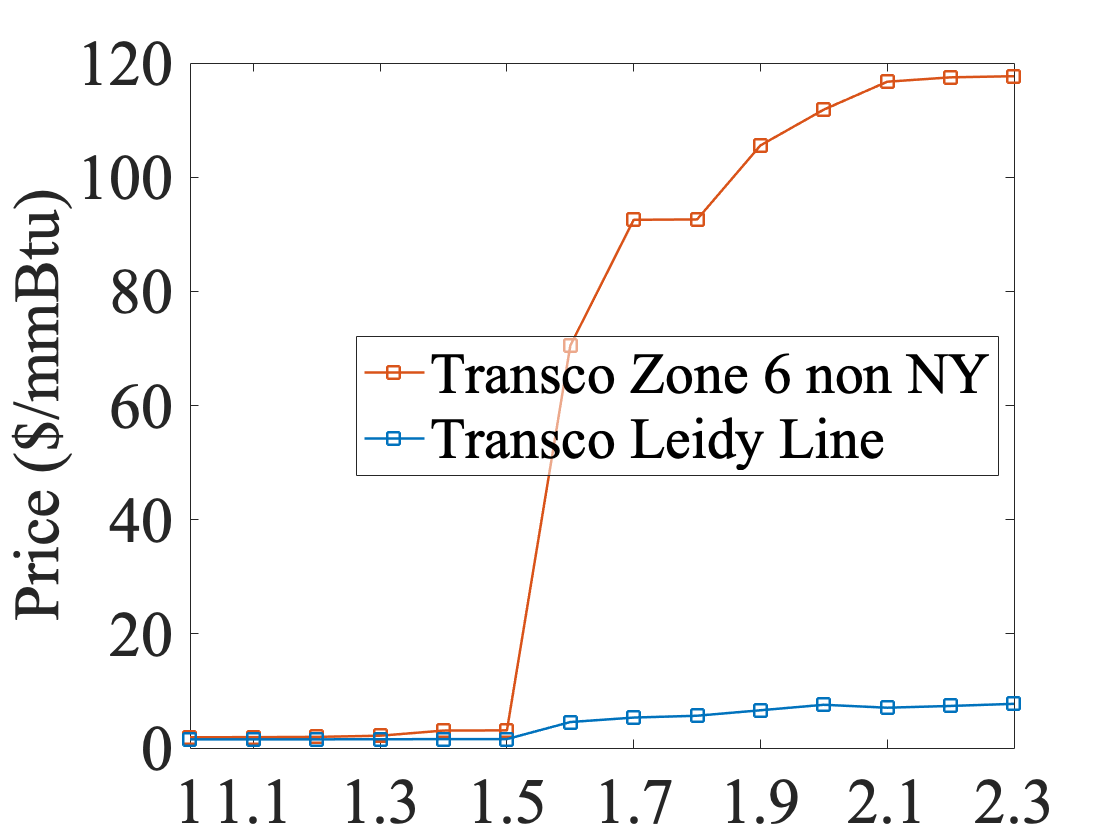}}
	\caption{Results for the Highly-Stressed Electrical Power System ($\eta_p = 1.6$).} \label{fig:pStress:17}
	\vspace{-0.25cm}
\end{figure}

The behaviors of (A) and (B) in the normal, stressed, and
highly-stressed power systems are compared in Figures
\ref{fig:pStress:1}, \ref{fig:pStress:13}, and \ref{fig:pStress:17}
respectively. In each figure, (a) and (c) display the system costs and
natural gas prices of (A), and (b) and (d) display those of (B). More
precisely, (a) and (b) present the total cost breakdown in terms of
the cost of electrical power system, the cost of the gas system, and
the economic loss from invalid bids. (c) and (d) depict the natural
gas zonal prices in each pricing zone.\footnote{Note that, as $\eta_g$
  increases, the total cost of (A) always increases, while the cost of
  (B) temporarily decreases sometimes. This is due to the presense of
  optimality gaps for some hard instances.}

Figures \ref{fig:pStress:1:cost:A} and \ref{fig:pStress:1:price:A}
show that the gas system cost gradually increases as $\eta_g$
increases up to $1.7$, then it grows rapidly from $\eta_g =1.8$ on. The rapid increase is due to load shedding (see Section
\ref{rem:price}) and leads to natural gas price spikes in Transco Zone
6 non NY. The large difference between the prices in Zone 6 and Leidy
Line indicates that the load shedding occurs due to the lack of
transmission capacity between these two points, not because of a lack
of gas supply. Due to the gas price spike in Transco Zone 6 non NY,
some bids of GFPPs become invalid and incur some losses, which
increases the total cost. On the other hand, for (B), the electrical
power system cost is slightly higher than for (A), but it does not
incur any economic loss from invalid bids and the overall cost is
lower.  Observe also that model (A) captures the same behavior as in
the 2014 polar vortex.  Additionally, observe that the gas price in
the Zone 6 region is also exhibiting sharp increases in model
(B). However, this peak 
has significantly less impact for (B) given the different commitment
decisions.

The differences in behavior between systems (A) and (B) become clearer
as the load increases in the electrical power system. For the stressed
power system, displayed in Figure \ref{fig:pStress:13}, the difference
between the total cost of (A) and (B) becomes very large: There are
many invalid bids for (A), which puts the reliability of the power
system at high risk and induces an electricity price peak. The price
of gas and the economic losses both increase significantly in (A) and
the increases start at stress level 1.5 for the gas network.  In
contrast, (B) maintains a reliable operation independently of the stress
imposed on the natural gas system. The price of gas increases
obviously but less than in (A) and the cost of the power system
remains stable. The peak in gas price only starts at stress level
1.7, showing that (B) delays the impact of congestion in the gas
networks by making better commitment decisions.

Figure \ref{fig:pStress:17} shows the benefits of (B) over (A) become
even more substantial when both systems are highly stressed. Observe
that the cost of the electrical power system remains stable once again
in (B) and that the cost of the gas network increases reasonably. In
contrast, Model (A) exhibits significant increases in gas prices and
economic cost from invalid bids. These results indicate that bringing
gas awareness in unit commitment brings significant benefits in
congested networks. By choosing commitment decisions that ensure bid
validity, the UCGNA brings substantial cost and reliability benefits
for congested situations like the 2014 polar vortex.

\begin{table}[!t]
\centering
	\caption{Statistics on Committed Generators for the Stressed Electrical Power System ($\eta_p = 1.6$): The first 7 columns display the number of committed generators with respect to its fuel type, where (O) Oil, (C) Coal, (G) Gas, (H) Hydro, (R) Refuse, (N) Nuclear, (E) Others, and the last two columns show the number of committed GFPPs in each pricing zone, where (T) Transco Zone 6 Non NY and (L) Transco Leidy Line. }\label{table:commit_gens_16}
	\begin{tabular}[h]{cccccccc|cc} 
			\hline
			$\eta_g$ & (O) & (C) & (G) & (H) & (R) & (N) & (E) & (T) & (L)\\
			\hline
			 1.0 & 7	& 6 &	12 &	11 &	0 &	12 &	3 & 8 & 4\\
			 1.6 & 8	& 6	& 10 &	11 &	0 &	13 &	3  & 6 & 4\\
			 2.3 & 9 &	6 &	9 &	11 &	0 &	13 &	3 & 4 & 4\\
			\hline
	\end{tabular}
\end{table}

\begin{figure*}[!t]
	\centering
	\includegraphics[width=0.55\textwidth]{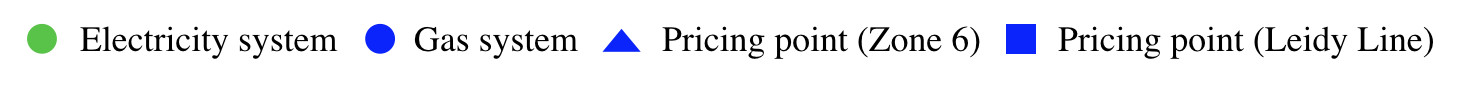}
	\subfloat[Number of committed GFPPs (A).\label{fig:UC:GFPP:A}] {\includegraphics[width=0.48\textwidth]{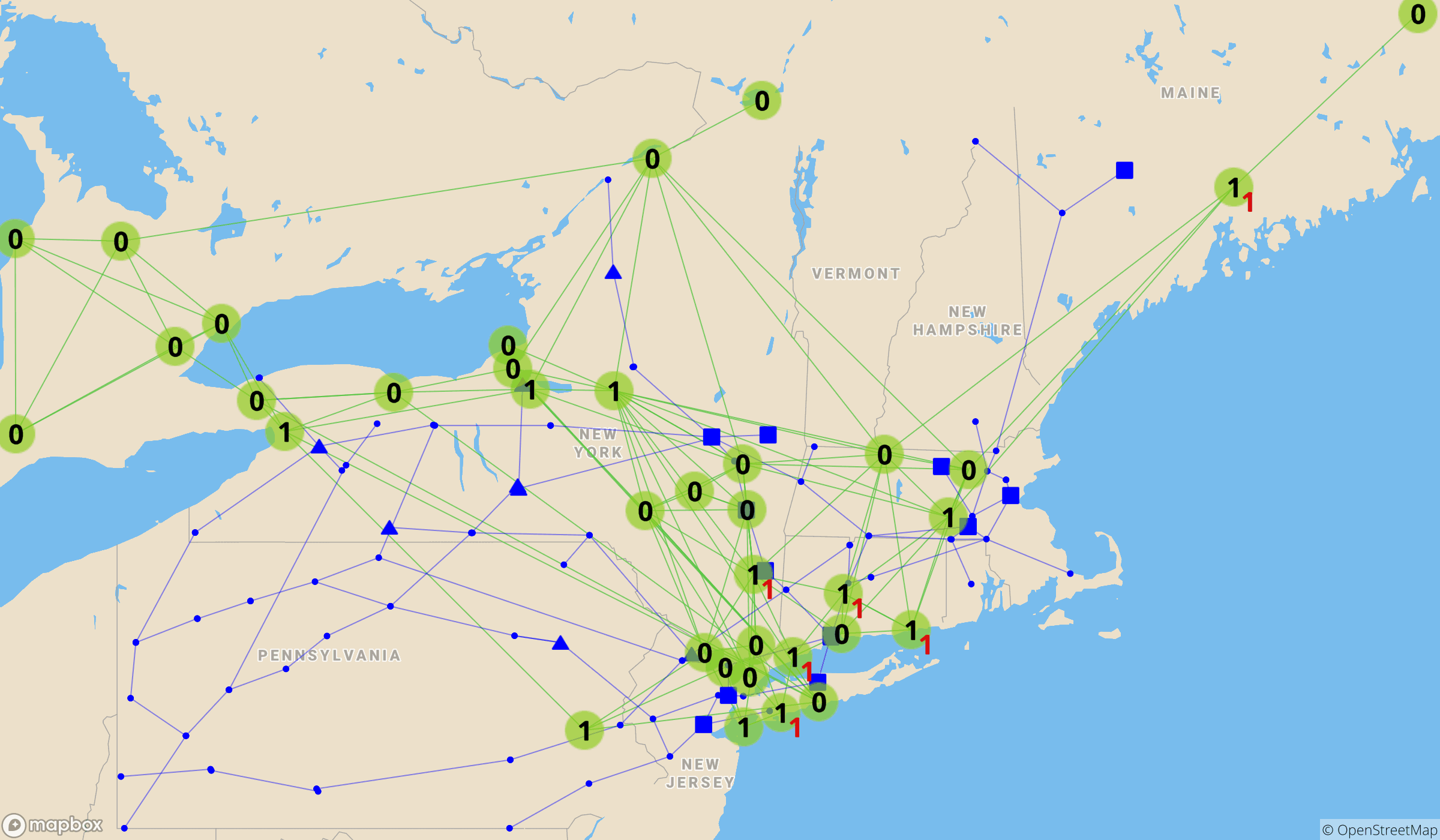}}\quad
	\subfloat[Number of committed non-GFPPs (A).\label{fig:UC:nGFPP:A}]  {\includegraphics[width=0.48\textwidth]{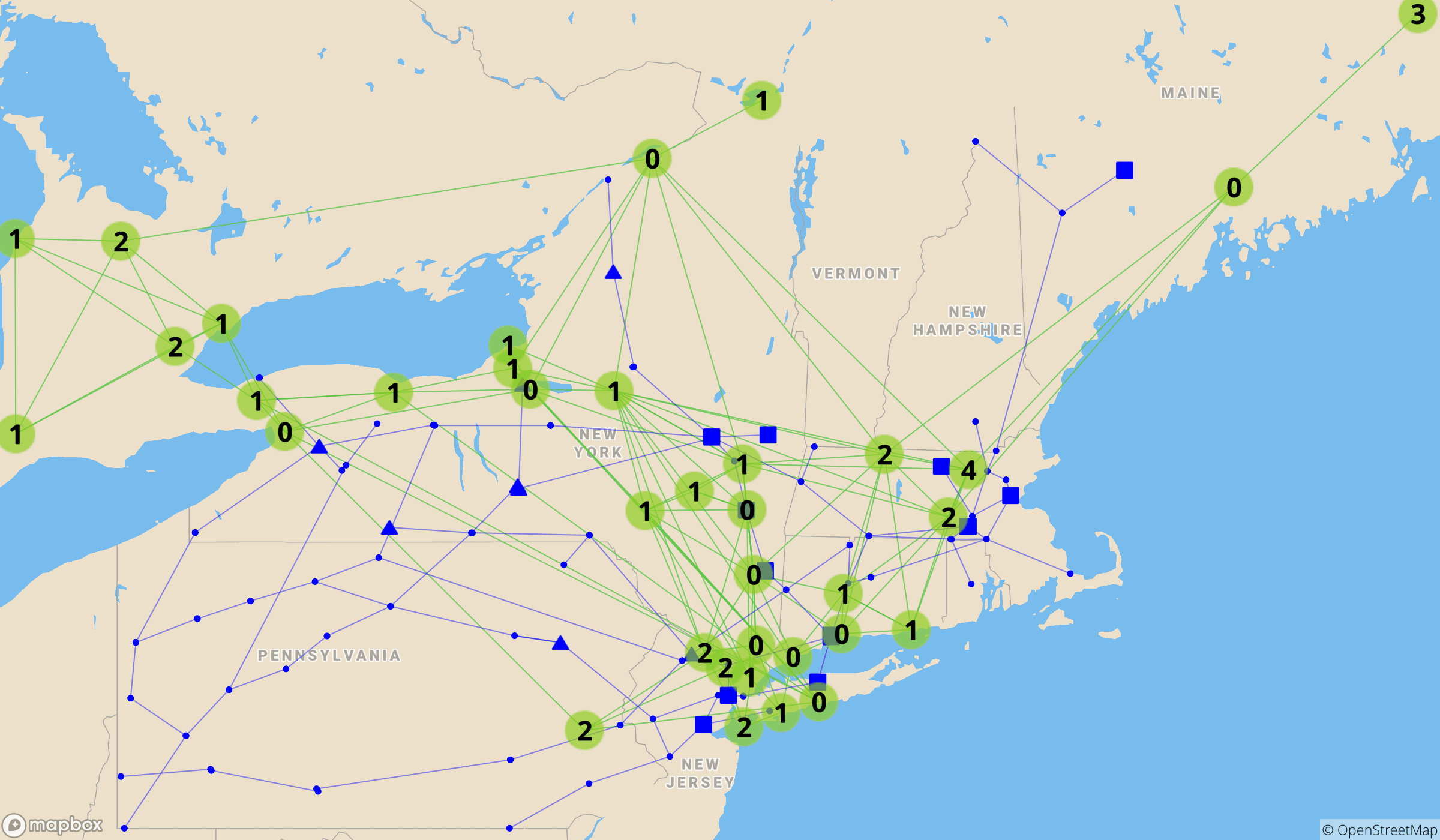}}\\
	\subfloat[Number of committed GFPPs (B).\label{fig:UC:GFPP:B}] {\includegraphics[width=0.48\textwidth]{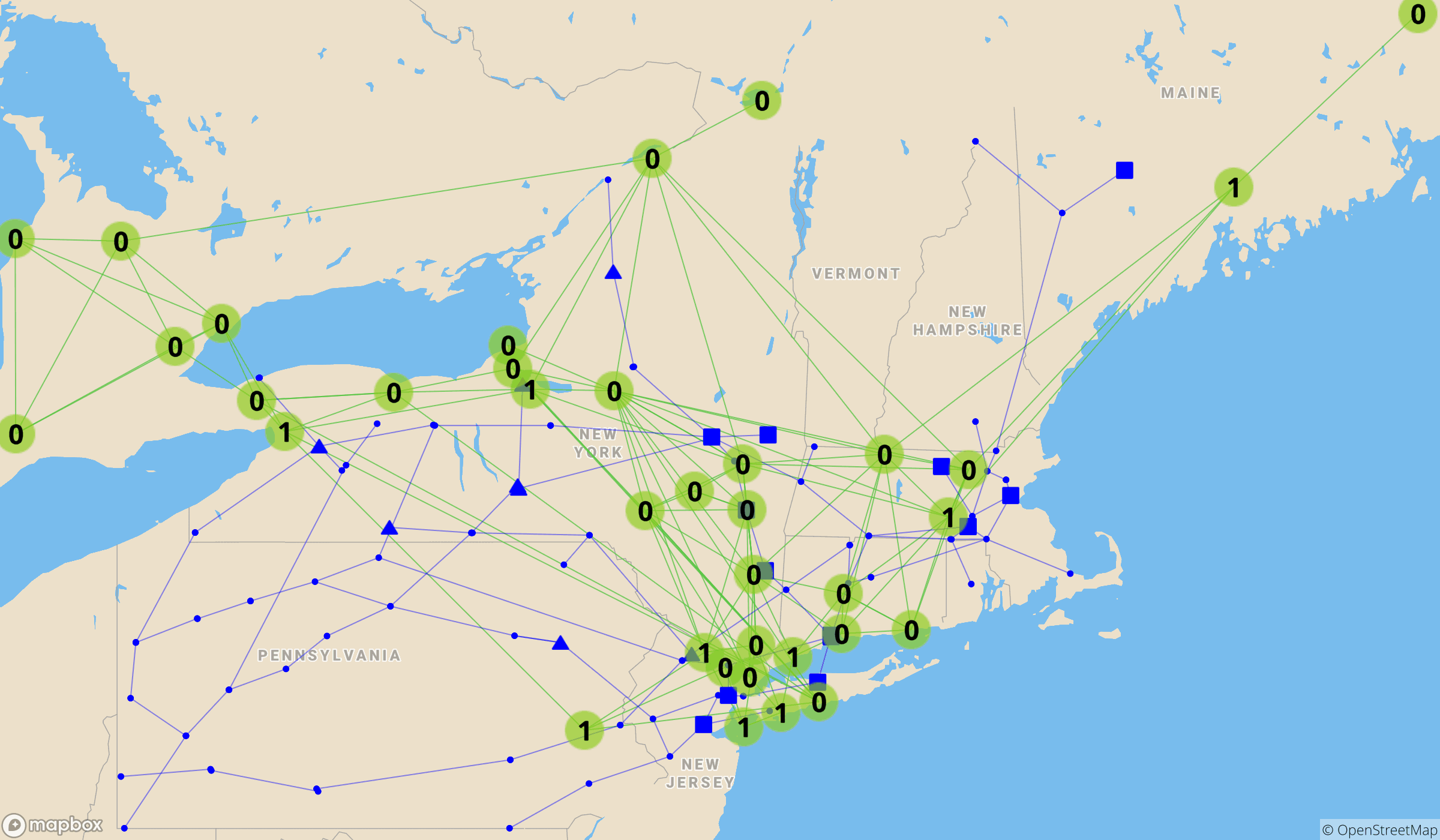}}\quad
	\subfloat[Number of committed non-GFPPs (B).\label{fig:UC:nGFPP:B}] {\includegraphics[width=0.48\textwidth]{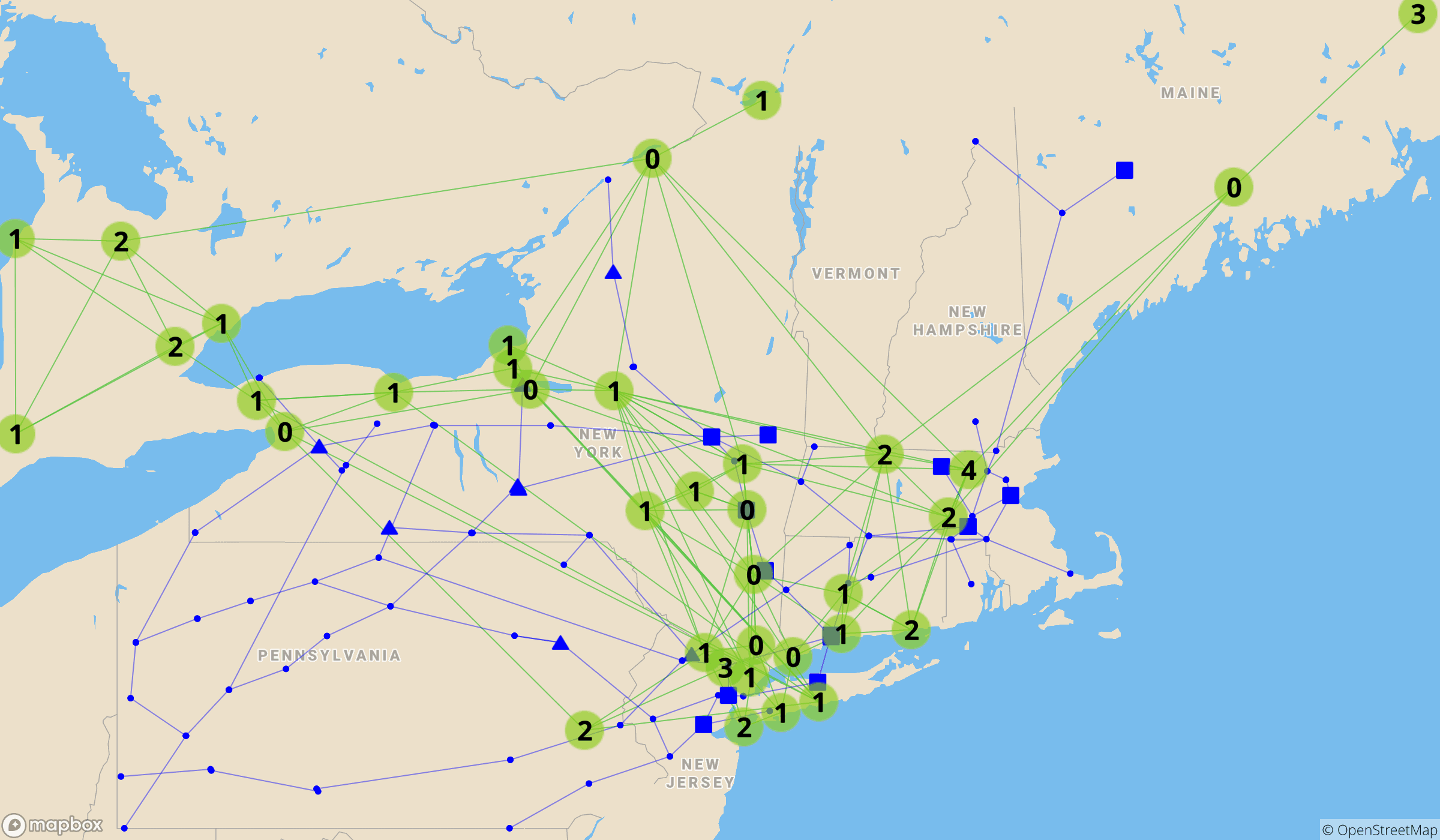}}\\
	\caption{Results for the Highly-Stressed Condition ($\eta_p,\eta_g$) = (1.6,2.3).}
	\label{fig:UC}
\end{figure*}

\begin{table}[!t]
	\caption{Solution Statistics for (B).}\label{table:obj}
	\begin{tabular}[h]{c|cc|cc|cc} 
			\hline
			\diagbox{$\eta_g$}{$\eta_p$} & \multicolumn{2}{c}{1} & \multicolumn{2}{c}{1.3} & \multicolumn{2}{c}{1.6}\\
			\cline{2-3} \cline{4-5} \cline{6-7}
			   & (i) & (ii) & (i) & (ii) & (i) & (ii)\\
			\hline
			  1	&	255301.0	&	0.0	&	332123.0	&	0.0	&	415315.0	&	0.0	\\
1.1	&	256502.0	&	0.0	&	333333.0	&	0.0	&	416530.0	&	0.0	\\
1.2	&	257706.0	&	0.0	&	334548.0	&	0.0	&	417759.0	&	0.0	\\
1.3	&	258915.0	&	0.0	&	335776.0	&	0.0	&	419015.0	&	0.0	\\
1.4	&	260132.0	&	0.0	&	337036.0	&	0.0	&	420548.0	&	0.0	\\
1.5	&	261364.0	&	0.0	&	338564.0	&	0.0	&	423466.0	&	0.0	\\
1.6	&	262613.0	&	0.0	&	342066.0	&	0.3	&	439254.0	&	2.1	\\
1.7	&	264019.0	&	0.0	&	361089.0	&	3.5	&	463746.0	&	2.0	\\
1.8	&	278679.0	&	1.8	&	379532.0	&	3.2	&	489011.0	&	6.2	\\
1.9	&	296251.0	&	1.3	&	408407.0	&	3.3	&	524533.0	&	7.4	\\
2	&	317619.0	&	0.0	&	430415.0	&	4.2	&	519026.0	&	3.7	\\
2.1	&	329801.0	&	0.0	&	460127.0	&	4.3	&	596449.0	&	5.0	\\
2.2	&	358828.0	&	0.0	&	497952.0	&	4.0	&	635128.0	&	5.0	\\
2.3	&	405022.0	&	0.0	&	537874.0	&	0.0	&	672876.0	&	0.0	\\
			\hline
	\end{tabular}
\end{table}

The great cost and reliability benefits of (B) are owing to better
commitment decisions that anticipate the future state of the gas
system. Table \ref{table:commit_gens_16} summarizes some statistics on
committed generators under the highly stressed power system. As the
gas load increases, some of the GFPPs in (T) are no longer committed
and the lost generation is replaced by generators of different types
or GFPPs with reasonable bid prices. More specifically, Figure
\ref{fig:UC} shows the commitment decision of (A) and (B) for
($\eta_p,\eta_g$) = (1.6,2.3). The numbers in black in Figures
\ref{fig:UC:GFPP:A} and \ref{fig:UC:GFPP:B} report the number of
committed GFPPs on the corresponding bus; Those in Figures
\ref{fig:UC:nGFPP:A} and \ref{fig:UC:nGFPP:B} display the number of
committed non-GFPPs. In Figure \ref{fig:UC:GFPP:A}, the numbers in red
on the bottom right corner of some buses represent the number of
committed GFPPs located at the bus without bid validity. Most invalid
GFPPs in Figure \ref{fig:UC:GFPP:A} are turned off in Figure
\ref{fig:UC:GFPP:B} and replaced by some non-GFPPs as Figure
\ref{fig:UC:nGFPP:B} indicates.

Finally, Table \ref{table:obj} summarizes the objective value and the
optimality gap of (B) for each instance. For 16 out of 42 instances,
the algorithm times out (wall-clock limit time of 1 hour) and it
reports sub-optimal solutions whose optimality gaps are presented in
columns denoted by (ii). It is important to stress however that even
sub-optimal solutions to the UCGNA bring significant benefits for
gas-grid networks as shown previously.

\section{Conclusion}\label{sec:conclusion}

The 2014 polar vortex showed how interdependencies between the
electrical power and gas networks may induce significant economic
and/or reliability risks under heavy congestion. This paper has
demonstrated that these risks can be effectively mitigated by making
unit commitment decisions informed by the physical and economic
couplings of the gas-grid network. The resulting Unit Commitment with
Gas Network Awareness (UCGNA) model builds upon the standard unit
commitment used in current practices but also reasons about the
feasibility of gas transmission feasibility and the profitability of
committed GFPPs. In particular, the UCGNA introduces bid-validity
constraints that ensure the economic viability of committed GFPPs,
whose marginal bid prices must be higher than the marginal natural gas
prices. The UCGNA is a three-level model whose bid validity
constraints operate on the dual variables of flux conservation
constraints in the gas network, which represent the marginal cost of
gas for producing a unit of electricity. It can be formulated as a
Mixed-Integer Second-Order Cone Program (MISOCP) and solved using a
dedicated Benders decomposition approach. The case study, based on a
modeling of the gas-grid network in the North-East of the United
States, shows that the UCGNA has significant benefits compared to the
existing operations: It is capable to ensure valid bids even at
highly-stressed levels, while only increasing the cost of gas and
electricity in a reasonable way. In contrast, the existing operating
practices induce significant economic losses and gas price increases.

Future research will be devoted to further improve the solution
techniques to solve the UCGNA and, in particular, the use of cut
bundling and Pareto-optimal cuts.

\section*{Acknowledgment}

This research was partly supported by an NSF CRISP Award (NSF-1638331)
``Computable Market and System Equilibrium Models for Coupled
Infrastructures''.

\ifCLASSOPTIONcaptionsoff
  \newpage
\fi

\newpage


%


\appendices
\section{Proof of Theorem \ref{theo:TL}}
\label{appendix:pf}

\begin{IEEEproof}
By strong duality of the third-level optimization in
Problem \eqref{prob:TL}, Problem \eqref{TL:lower} is equivalent to:

\begin{subequations}\fontsize{9}{9}\selectfont
\begin{alignat}{3}
(\x_p, \y_g) = &  \underset{\x_p \ge 0, \y_g }{\mbox{argmin}} && c_p^T \x_p\\ 
& \mbox{s.t.} &&  A \x_p + B \z_p \ge b, \label{TL:lower2:upper}\\
& &&	{\begin{array}{lcl}
  \y_g = &\underset{\x_g \in \K, \y_g \ge 0}{\mbox{argmin}} \ & c_g^T \x_g \\ 
& \mbox{s.t.}  & D_p \x_p + D_g \x_g \ge d,\\
& & \y_g^T (d - D_p \x_p) \ge c_g^T \x_g, \\
& &  \y_g^T D_g \preceq_{\K^*} c_g. \end{array}} \label{TL:lower2:lower}
\end{alignat}
\label{prob:TL:lower2}
\end{subequations}
where $\K^*$ denotes the dual cone of $\K$. The first and third
constraints of Problem \eqref{TL:lower2:lower} state the primal and
dual feasibility of the third-level problem, while the second
constraint ensures their optimality.
	

Equation \eqref{TL:lower2:upper} (i.e., the constraint of the upper
level problem of Problem \eqref{prob:TL:lower2}) does not involve the
lower-level variables (i.e., $\x_g$ and $\y_g$ of
Problem \eqref{TL:lower2:lower}), which means the upper-level solution
is not affected by the solutions to the lower-level problem.
Problem \eqref{prob:TL:lower2} can thus be solved in two steps: (i)
solve the upper-level problem and obtain $\bar x_p$, (ii) solve the
lower-level problem with $\x_p$ fixed as $\bar x_p$ and obtain $\bar
y_g$. Accordingly, Problem \eqref{prob:TL:lower2} can be expressed
with a Lexicographic function as follows:
	
\begin{subequations}\fontsize{9}{9}\selectfont
			\begin{alignat}{3}
			(\x_p, \y_g) = &  \underset{\x_p \ge 0, \x_g \in \K, \y_g \ge 0}{\mbox{argmin}} && <c_p^T \x_p, c_g^T \x_g> \\ 
			&\mbox{s.t.} && A \x_p + B \z_p \ge b, \\
			& &&	D_p \x_p + D_g \x_g \ge d, \\
			& && \y_g^T (d - D_p \x_p) \ge c_g^T \x_g,\\
			& && \y_g^T D_g \preceq_{\K^*} c_g.
			\end{alignat}
			\label{prob:Lex}
\end{subequations}

\noindent
The optimal solution $(\bar x_p, \bar x_g, \bar y_g)$ of
Problem \eqref{prob:Lex} satisfies the following conditions:
	
\begin{subequations}\fontsize{9}{9}\selectfont
\begin{alignat}{3}
\bar x_p =\underset{\x_p \ge 0,\x_g \in \K}{\mbox{argmin}} \ &    c_p^T \x_p \\
\mbox{s.t. }&    A \x_p \ge b - B \z_p, \\
	    &	 D_p \x_p + D_g \x_g \ge d.
\end{alignat} 
\label{prob:Lex:first}
\end{subequations}
\begin{subequations}\fontsize{9}{9}\selectfont
\begin{alignat}{3}
(\bar x_g, \bar y_g) = \underset{\x_g \in \K, \y_g \ge 0}{\mbox{argmin}} \ & c_g^T \x_g \\ 
\mbox{s.t.} \ & D_g \x_g \ge d - D_p \bar x_p, \\
& \y_g^T (d - D_p \bar x_p) \ge c_g^T \x_g,\label{Lex:2nd:st}\\
& \y_g^T D_g \preceq_{\K^*} c_g.
\end{alignat} 
\label{prob:Lex:2nd}
\end{subequations}
	
\noindent
Observe that any feasible $(\hat x_g, \hat y_g)$ of
Problem \eqref{prob:Lex:2nd} is optimal. By strong duality, ($\hat
x_g, \hat y_g$) satisfies the following conditions:
		
\begin{subequations}\fontsize{9}{9}\selectfont
\begin{alignat}{3}
\hat x_g = \underset{\x_g \in \K}{\mbox{argmin}} \ & c_g^T \x_g\\
\mbox{s.t. }&  D_g \x_g \ge d - D_p \bar x_p.
\end{alignat} 
\label{prob:Lex:2nd:primal}
\end{subequations}
\begin{subequations}\fontsize{9}{9}\selectfont
\begin{alignat}{3}
\hat y_g = \underset{\y_g \ge 0}{\mbox{argmax}} \ &  \y_g^T (d - D_p \bar x_p)\\
\mbox{s.t. }& \y_g^T D_g\preceq_{\K^*} c_g.
\end{alignat} 
\label{prob:Lex:2nd:dual}
\end{subequations}
\noindent		
Since Problem \eqref{prob:Lex:2nd:primal} is a relaxation of
Problem \eqref{prob:Lex:2nd} and $\hat x_g$, paired with $\hat y_g$,
is feasible for Problem \eqref{prob:Lex:2nd}, $(\hat x_g, \hat y_g)$ is
optimal to Problem \eqref{prob:Lex:2nd}. As a result, for $\alpha \in
(0,1)$, Problem \eqref{prob:TL} can be approximated by

\begin{subequations}\fontsize{9}{9}\selectfont
\begin{alignat}{3}
\min \ & \alpha h^T \z_p + \alpha c_p^T \x_p  + (1-\alpha) c_g^T \x_g  \label{Lex_alpha:obj:leader}\\
\mbox{s.t.} \ & \z_p \in \mathcal Z, \\
& (\x_p, \x_g, \y_g) = \mbox{\fontsize{9}{9}\selectfont Primal \& dual opt. sol. of } \eqref{prob:Lex_alpha_lower}, \label{Lex_alpha:lower}\\
& \frac{1}{1-\alpha} E \y_g  + M \z_p \ge h, \label{Lex_alpha:leader},\\
&\x_p \ge 0, \x_g \in \K, \y_p \ge 0, \y_g \ge 0, \\
&\z_p \in \{0,1\}^m.
\end{alignat}
\label{prob:Lex_alpha}
\end{subequations}
\noindent
where the low-level problem in Equation \eqref{Lex_alpha:lower} is
\begin{subequations}\fontsize{9}{9}\selectfont
\begin{alignat}{3}
\underset{\x_p \ge 0, \x_g \in \K}{\min} &  \alpha c_p^T \x_p  + (1-\alpha)  c_g^T \x_g\\
\mbox{s.t. }&    A \x_p  + B z_p\ge b,\\ 
&  D_p \x_p + D_g \x_g \ge d. \label{eq:Lex:lower:2}
\end{alignat}
\label{prob:Lex_alpha_lower}
\end{subequations}
	
Problem \eqref{prob:Lex_alpha_lower} is an approximation of
Problem \eqref{prob:Lex}, where $\y_g$ is obtained by the dual
solution associated with Equation \eqref{eq:Lex:lower:2}. Hence, by
strong duality of Problem \eqref{prob:Lex_alpha_lower},
Problem \eqref{prob:our} is equivalent to
Problem \eqref{prob:Lex_alpha}.

It remains to show that Problem \eqref{prob:our} is indeed an
asymptotic approximation of Problem \eqref{prob:TL}.  Replacing $y_p$
with $y_p / \alpha$ and $y_g$ with $y_g / (1-\alpha)$ in
Problem \eqref{prob:our} gives the following equivalent problem:
	
\begin{subequations}\fontsize{9}{9}\selectfont
\begin{alignat}{3}
\min \ & \alpha h^T \z_p + \alpha c_p^T \x_p  +  (1-\alpha) c_g^T \x_g \label{our2:obj}\\
\mbox{s.t.} \ & 		\z_p \in \mathcal Z, \\
&  	A \x_p + B \z_p \ge b,\label{our2:1st:primal}\\
& D_p \x_p + D_g \x_g \ge d,\label{our2:follower}\\
&  \y_p^T (b - Bz_p) - c_p^T \x_p \ge \frac{1-\alpha}{\alpha}\left[c_g^T \x_g  - \y_g^T d \right], \label{our2:strong}\\
& \y_g^T D_g  \preceq_{\K^*} c_g^T ,\label{our2:2nd:dual_geas}\\
& \y_p^T A + \frac{1-\alpha}{\alpha} \y_g^T  D_p  \le  c^T_p,\label{our2:1st:dual_geas}\\
& E \y_g + M \z_p \ge h,\\
&\x_p \ge 0, \x_g \in \K, \y_p \ge 0, \y_g \ge 0, \\
&\z_p \in \{0,1\}^m.
\end{alignat}
\label{prob:our2}
\end{subequations}
	
Let $P(\hat z_p)$ and $\widehat P(\hat z_p)$ denote
Problems \eqref{prob:Lex} and \eqref{prob:our2} in which the
binary variables $\z_p$ are fixed to some $\hat
z_p \in \{0,1\}^m$.  Let $(\hat x_p, \hat x_g, \hat y_p, \hat
y_g)$ be the optimal solution of $\widehat P(\hat z_p)$. Note
that, as $\alpha \rightarrow 1$, Equations \eqref{our2:strong}
and \eqref{our2:1st:dual_geas} become as follows:
	
\begin{subequations}\fontsize{9}{9}\selectfont
\begin{alignat}{3}
&  \y_p^T (b - B \hat z_p) \ge  c_p^T \x_p, \label{1st:inf_alpha:strong}\\
&  \y_p^T A \le  c^T_p,\label{1st:inf_alpha:dual_geas}
\end{alignat}
\label{prob:1st:inf_alpha}
\end{subequations}

\noindent
which implies that $\hat y_p$ and $\hat x_p$ approximate the optimal
primal and dual solutions of Problem \eqref{prob:Lex:first} when $\z_p$
is fixed as $\hat z_p$. This is because $\hat x_p$ is feasible
for \eqref{prob:Lex:first} (by Equation \eqref{our2:1st:primal}), $\hat
y_p$ becomes feasible to the dual of Problem \eqref{prob:Lex:first} as
$\alpha$ approaches 1 (by Equation \eqref{1st:inf_alpha:dual_geas}),
and together they satisfy the strong duality condition of
Equation \eqref{1st:inf_alpha:strong} as $\alpha$ becomes closer to 1
(by Equation \eqref{1st:inf_alpha:strong}). Therefore, as
$\alpha \rightarrow 1$, $(\hat x_p, \hat y_p)$ becomes a feasible
solutions of $P(\hat z_p)$ and has the same optimal objective value.

Moreover, combining Equations \eqref{our2:strong} and \eqref{our2:1st:dual_geas}
gives 
	
\begin{subequations}\fontsize{9}{9}\selectfont
\begin{alignat}{4}
& &&( \mbox{Equation } \eqref{our2:strong}) -  \hat x_p  \times (\mbox{Equation }  \eqref{our2:1st:dual_geas})\nonumber\\
& \rightarrow \quad&&  \hat y_p^T  (b- B \hat z_p - A \hat x_p) +  \frac{1-\alpha}{\alpha}\hat y_g^T  (d - D_p \hat x_p) \ge    \frac{1-\alpha}{\alpha} c_g^T \hat x_g\nonumber\\
& \rightarrow &&  \hat y_g^T  (d - D_p \hat x_p) \ge   c_g^T \hat x_g,\label{2nd:strong}
\end{alignat}		
\end{subequations}
\noindent
where the last derivation follows from
Equation \eqref{our2:1st:primal} and $y_g \ge 0$. Therefore, $\hat
x_g$ and $\hat y_g$ are the optimal solutions of
Problem \eqref{prob:Lex:2nd} when $x_p$ is fixed as $\hat x_p$ (since
its feasibility is guaranteed by Equations \eqref{our2:follower}
and \eqref{our2:2nd:dual_geas}, while the optimality is guaranteed by
Equation \eqref{2nd:strong}).

In summary, $\hat x_p$ is an approximate solution of $P(\hat z_p)$
that becomes increasingly close to the optimal solution of
Problem $P(\hat z_p)$ as $\alpha \rightarrow 1$, and $\hat
y_g$ is the exact response of the follower with respect to
$\hat x_p$ for any $\alpha \in (0,1)$. Therefore, the
approximation may sacrifice the leader's optimality when
$\alpha$ is not large enough, but it always gives a feasible
solution.
\end{IEEEproof}

\end{document}